\documentclass[11pt]{article}

\usepackage{a4wide}
\usepackage{amstext}
\usepackage{amssymb}
\usepackage{latexsym}
\newcommand\RR{{\mathbb R}}

\newcommand\QQ{{\mathbb Q}}
\newcommand\NN{{\mathbb N}}
\newcommand\ZZ{{\mathbb Z}}

\newcommand\HH{{\mathbb H}}
\def\d={\,:=\,}

\font\frakten=eufm10
\newfam\frakfam
\textfont\frakfam=\frakten

\newtheorem{thm}{Theorem}[section]
\newtheorem{lemma}[thm]{Lemma}
\newtheorem{cor}[thm]{Corollary}
\newtheorem{prop}[thm]{Proposition}
\newtheorem{Defn}[thm]{Definition}
\newtheorem{Ex}[thm]{Example}
\newtheorem{Rem}[thm]{Remark}
\newtheorem{Exs}[thm]{Examples}
\newtheorem{Rems}[thm]{Remarks}
\newtheorem{Defrem}[thm]{Definition and Remark}
\newtheorem{Remnt}[thm]{}
\newenvironment{defn}
 {\begin{Defn} \begin{rm}} {\end{rm} \hfill $\Box$ \end{Defn}}

\newenvironment{ex}
 {\begin{Ex} \begin{rm}} {\end{rm} \hfill $\Box$ \end{Ex}}

\newenvironment{rem}
 {\begin{Rem} \begin{rm}} {\end{rm} \hfill $\Box$ \end{Rem}}
\newenvironment{rems}
 {\begin{Rems} \begin{rm}} {\end{rm} \hfill $\Box$ \end{Rems}}

\newenvironment{prf} {{\bf Proof.}}{\hfill $\Box$}

\def\be{\begin{equation}}
\def\ee{\end{equation}}
\def\bea{\begin{eqnarray}}
\def\eea{\end{eqnarray}}
\def\text#1{{\rm #1}}

\def\lb{\hbox{$\lambda$\hskip-5pt$\raise2pt\hbox{\char'040}$\hskip1.5pt}}

\begin{document}

\title{Sampling theorems for the Heisenberg group}
\author{Hartmut F\"uhr \\ GSF-IBB\\ 85764 Neuherberg\\ Germany}
\date{\today}
\maketitle

%--------------------------------------------------------------

\begin{abstract}
 In the first part of the paper a general notion of sampling expansions
 for locally compact groups is introduced, and its close relationship
 to the discretisation problem for generalised wavelet transforms is
 established. In the second part, attention is focussed on the
 simply connected nilpotent Heisenberg group $\HH$. We derive
 criteria for the existence of discretisations and sampling expansions
 associated to lattices in $\HH$. 
 Analogies and differences to the sampling theorem over the reals
 are discussed, in particular a notion of bandwidth on $\HH$
 will figure prominently. The main tools for the characterisation are
 the Plancherel formula of $\HH$ and the theory of Weyl-Heisenberg frames.
 In the last section we compute an explicit example.
\end{abstract}

\section{Introduction}\label{sec:Intro}

Let us first introduce some notation and definitions.
$G$ denotes a locally compact topological group, with left Haar
measure $\mu_G$, and associated ${\rm L}^2$-space ${\rm L}^2(G)$.
The {\bf left regular representation} $\lambda_G$ acts on ${\rm L}^2(G)$
by left translation. For a function $g$ on $G$, $g^*$ is defined
as $g^*(x) = \overline{g(x^{-1})}$. Representations of $G$ are understood
to be strongly continuous and unitary.

In this paper we wish to discuss the existence of sampling theorems
analogous to the famous Whittaker-Shannon-Kotel´nikov Theorem.
The setting we study is group-theoretic: Given a discrete subset 
$\Gamma \subset G$,
we wish to reconstruct certain functions $g \in {\rm L}^2(G)$ from
their values sampled at $\Gamma$. We formalise this in the following 
notion.

\begin{defn}
\label{def:samp_space}
 Let $G$ be a locally compact group, $\Gamma \subset G$. Let
 ${\cal H} \subset {\rm L}^2(G)$ be a leftinvariant
 closed subspace of ${\rm L}^2(G)$ consisting of continuous functions. 
 We call ${\cal H}$ a {\bf sampling space (with respect to $\Gamma$)}
 if it has the following two properties:
 \begin{enumerate}
 \item[(i)] There exists a constant $c_{\cal H}>0$, such that for all
 $f \in {\cal H}$,
 \[ \sum_{\gamma \in \Gamma} |f(\gamma) |^2 = c_{\cal H} \| f \|_2^2 ~~.\]
 In other words, the restriction mapping 
 $R_{\Gamma}: {\cal H} \ni f \mapsto (f |_{\Gamma}) \in
 \ell^2(\Gamma)$ is a scalar multiple of an isometry.
 \item[(ii)] There exists $S \in {\cal H}$ such that every $f \in {\cal H}$
 has the expansion
 \begin{equation} \label{eqn:samp_exp}
 f(x) = \sum_{\gamma \in \Gamma} f(\gamma) S(\gamma^{-1} x)~~, 
 \end{equation}
 with convergence both in ${\rm L}^2$ and uniformly. 
 \end{enumerate}
 The function $S$ from condition (ii) is called {\bf sinc-type function}.
\end{defn}

The definition is modelled after the following, prominent example:
\begin{ex} [{\bf Whittaker, Shannon, Kotel'nikov}]
Let $G = \RR$, $\Gamma = \ZZ$ and 
 \[ {\cal H} = \{ f \in {\rm L}^2(\RR) : {\rm supp}(\widehat{f}) \subset
 [-0.5,0.5]  \}~~. \]
Then ${\cal H}$ is a sampling subspace, with 
\[ S(x) = {\rm sinc}(x) = \frac{\sin (\pi x)}{\pi x} ~~.\]
\end{ex}

One of the aspects of the Whittaker-Shannon-Kotel'nikov
sampling theorem has not been covered by Definition \ref{def:samp_space}: 
The sampling expansion 
(\ref{eqn:samp_exp}) can be read as interpolating between 
values of $g$ at $\Gamma$. The original sampling theorem allows to interpolate
arbitrary $\ell^2$-sequences by 
elements of ${\cal H}$, i.e., the restriction map is {\em onto}.
We do not require this for sampling subspaces. It will become apparent in
the second part of the paper that the Heisenberg group allows a variety of 
sampling spaces,
but none with arbitrary interpolation. 

Sampling spaces will be discussed in connection with notions coming from
wavelet theory, such as admissible vectors and tight frames.
A system $(\eta_i)_{i \in I}$ of vectors in a Hilbert
space ${\cal H}$ is a {\bf tight frame}
if the {\bf coefficient operator} 
\begin{equation} \label{eqn:coeff_operator}
 T: {\cal H} \ni \varphi \mapsto (\langle \varphi, \eta_i\rangle )_{i \in I} 
 \in \ell^2(I) 
\end{equation}
is isometric up to a constant, or equivalently, if
\[ \varphi = \frac1c \sum_{i \in I} \langle \varphi, \eta_i \rangle
 \eta_i ~~,\]
holds for every $\varphi \in {\cal H}$. The sum converges (unconditionally) 
in the norm. A tight frame is called {\bf normalised} if the 
constant $c$ equals 1. It is clear that a tight frame is total in ${\cal H}$.
The following proposition collects some basic
used facts about tight frames. The facts and their proofs are widely known
and just repeated for the sake of completeness.
\begin{prop} \label{prop:basics_frames}
 Let $(\eta_i)_{i \in I} \subset {\cal H}$ be a tight frame
 with frame constant $c$.
 \begin{enumerate}
 \item[(a)] If ${\cal H}' \subset {\cal H}$ is a closed subspace and
 $P: {\cal H} \to {\cal H}'$ is the projection onto ${\cal H}'$, then
 $(P \eta_i )_{i \in I}$ is a tight frame of ${\cal H}'$ with frame
 constant $c$.
 \item[(b)] Suppose that $c=1$. Then $(\eta_i)_{i \in I}$ is an orthonormal
 basis iff $\| \eta_i \| = 1$, for all $i \in I$. 
 \item[(c)] If $\| \eta_i \| = \| \eta_j \|$, for all $i, j \in I$, then
 $\| \eta_i \|^2 \le c$.
 \item[(d)] For an arbitrary system $(\phi_i)_{i \in I} \subset {\cal H}$
 define the coefficient operator $T$ analogously to (\ref{eqn:coeff_operator}),
 with maximal possible domain, i.e. ${\rm dom} (T) = \{ z \in {\cal H} :
 (\langle z,\phi_i \rangle)_{i \in I} \in \ell^2(I) \}$ (which may
 be trivial). Then $T$ is a closed operator. 
 \item[(e)] $(\eta_i)_{i \in I}$ is an orthornormalbasis iff $c=1$ and the
 coefficient operator is onto.
 \end{enumerate}
\end{prop}

\begin{prf}
 Part $(a)$ follows from the fact that on ${\cal H}'$ the coefficient map
 associated to $(P \eta_i )_{i \in I}$ coincides with the coefficient map
 associated to $(\eta_i )_{i \in I}$. The ``only-if''-part
 of $(b)$ is clear. The ``if''-part follows from 
 \[ 1 = \| \eta_i \|^2 = \sum_{i \in I} | \langle \eta_i, \eta_j \rangle|^2
 = 1 + \sum_{i \not= j }  | \langle \eta_i, \eta_j \rangle|^2~~,
 \]
 whence $ \langle \eta_i, \eta_j \rangle$ vanishes for $i \not= j$.
 Part $(c)$ follows from a similar argument. The proof of part $(d)$
 is a straightforward application of the Cauchy-Schwarz inequality.
 The ``only if'' part of $(e)$ is obvious. For the converse let
 $\delta_i \in \ell^2(I)$ be the Kronecker delta at $i$. Then
 $\langle T^* \delta_i, \varphi \rangle = \langle \delta_i, T \varphi \rangle 
 = \langle \eta_i , \varphi \rangle$ for all $\varphi \in H$ implies
 $T^* \delta_i = \eta_i$, or $T \eta_i = \delta_i$ ($T$ is by
 assumption unitary), which is the
 desired orthonormality relation.
\end{prf}

In this paper, all frames of interest will be of the form
\[ \pi(\Gamma) \eta = ( \pi(\gamma) \eta )_{\gamma \in \Gamma} ~~,\]
where $\pi$ is a representation of $G$ (usually the restriction
of $\lambda_G$ to some leftinvariant closed subspace ${\cal H}
\subset {\rm L}^2(G)$), and $\Gamma$
a suitable discrete subset. The associated coefficient map is usually
called ``discrete wavelet transform''.

% However, as the case
%of the Heisenberg group shows, the additional freedom is necessary,
%since for the sampling sets under consideration, no sinc-type
%function exists which generates an ONB, but there do exist sinc-type
%functions which generate a normalised tight frame.

Let us give a short outline of the paper: We first establish a close
connection between the discretisation of (generalised) wavelet transforms
and sampling spaces, for general locally compact groups.
The discretisation problem for wavelet transforms and the
construction of sampling spaces are essentially equivalent. 
We then address these problems for a concrete example, the three-dimensional 
Heisenberg group $\HH$. Section 3 contains the basic results on
the Heisenberg group which are used in the following. In Section 4 we
formulate and discuss our main results. Our discussion aims at exposing 
analogies
as well as differences to the sampling theorem over the reals.
In particular, a notion of bandwidth will turn out to play
a prominent role. Sections 5 through 7 contain the proofs
of the main results. The main technical devices used for the
proofs are the Plancherel formula for the Heisenberg group and
the theory of Weyl-Heisenberg frames. In Section 8 we illustrate our techniques
by explicitly computing a sinc-type function on $\HH$.
An appendix contains notations and results used in connection
with the Plancherel formula, concerning Hilbert-Schmidt
operators, direct integrals and the like.

 For locally compact abelian groups, the paper by Kluv\'anek \cite{Kl} 
is considered to contain the definitive form of the sampling theorem.
In chapter 10 of \cite{HiSt}, which
covers these results, extensions to nonabelian groups
are mentioned as a natural but difficult question. To our 
knowledge, the paper by Dooley \cite{Do} is the only source
dealing with sampling expansions for nonabelian groups.
As was observed in \cite[chapter 10]{HiSt}, a straightforward generalisation 
of Kluv\'anek's approach is not possible,
since the dual lattice, which is the central device in the abelian case,
is not available in the context of nonabelian groups (the dual is
not a group). However Dooley's results indicated that the
Plancherel formula of a nonabelian group could be used to
study sampling expansions, and our paper provides further evidence.

The initial purpose of this paper was to investigate whether
the Plancherel formula could be used for the discretisation of 
continuous wavelet transforms. It was shown recently in \cite{AlFuKr,FuMa,Fu}
that Plancherel theory provides a natural framework for a unified treatment
of continuous wavelet transforms and Wigner functions, and
extending this framework to discretisation problems seems to be
an attractive project. Moreover, the similarity
between the discretisation problem and the Shannon sampling theorem has been
observed, more or less explicitly, by various authors, most notably by 
Feichtinger and Gr\"ochenig in their series of papers (see the bibliography
of \cite{Gr}). The approach via the Plancherel formula allows to address 
this analogy explicitly. We chose the Heisenberg group mainly as a test
case which should provide some orientation for more general settings.

\section{Discretised wavelet transforms and sampling
 expansions}
\label{sect:dwt_ss}

In this section we establish a close connection between sampling
expansions and the problem of discretising generalised wavelet
transforms. 
\begin{defn} Let $(\pi,{\cal H}_{\pi})$ be a representation of $G$.
$\eta \in {\cal H}_{\pi}$ is called {\bf admissible} if the
coefficient mapping
\[ V_{\eta} : \varphi \mapsto V_{\eta} \varphi ~~,\]
where
\[ V_{\eta} (\varphi) (x) = \langle \varphi, \pi(x) \eta \rangle ~~,\]
is an isometry ${\cal H}_{\pi} \to {\rm L}^2(G)$. The operator $V_{\eta}$ is
then called {\bf (generalized) continuous wavelet transform}, it
intertwines $\lambda_G$ with $\pi$. If an admissible vector exists,
$\pi$ is called {\bf square-integrable}. The admissibility 
condition implies the following reconstruction formula (to be read in
the weak sense)
\[ \phi = \int_G V_{\eta} \phi(x) ~ \pi(x) \eta ~
 d\mu_G(x)~~.
\]
\end{defn}
Note that there exist various definitions of square-integrable representations
in the literature. Perhaps the most common one is that $\pi$ is in the
{\em discrete series}, i.e., $\pi$ is irreducible and has
admissible vectors. Square-integrable representations, as defined here,
are (equivalent to) subrepresentations of the regular representation;
the converse does not hold in general, see \cite{Fu}.
If $\pi$ is the restriction of $\lambda_G$ to a leftinvariant
subspace, then $V_{\eta} \phi = \phi \ast \eta^*$. 
For groups with a type I regular representation there exist precise 
criteria for square-integrability which employ the Plancherel measure 
of such groups, see \cite{Fu}. For the Heisenberg group they are 
formulated in Theorem \ref{thm:adm_vec}.

The connection between generalised wavelet transforms and sampling
spaces is realised via the space $V_{\eta}({\cal H}_{\pi})$. The space
is characterised by the {\bf reproducing kernel relation}
\[ f = f \ast V_{\eta} \eta ~~,\]
which follows by elementary computation from the isometry property
of $V_{\eta}$. The convolution kernel $S = V_{\eta} \eta$ is a
{\bf (right) selfadjoint convolution idempotent} in ${\rm L}^2(G)$, i.e.
it fulfills $S = S \ast S = S^*$. Considering subspaces of the
form ${\cal H} = V_{\eta}({\cal H}_\pi)$, for admissible $\eta$,
is equivalent to studying selfadjoint convolution idempotents
in ${\rm L}^2(G)$. Indeed, given such an idempotent $S$, it
is straigthforward to check that $f \mapsto f \ast S$ is the orthogonal
projection onto a closed, leftinvariant subspace ${\cal H} = {\rm L}^2(G) 
\ast S$. With suitable identifications this space is easily seen
to be the image of a generalised wavelet transform: Simply pick
$S = S \ast S \in {\cal H}$ as admissible vector for the restriction
of $\lambda_G$ to ${\cal H}$, then the generalised wavelet transform
$V_S$ is the inclusion map.

The discretisation problem for a generalised wavelet transform
can be phrased as follows: Find a suitable (admissible)
$\eta$ and a sampling subset $\Gamma \subset G$, such that 
$\pi(\Gamma) \eta$ is a tight frame. Now the definition of a sampling
space implies that it is enough to ensure that $V_{\eta} ({\cal H}_{\pi})$
is a sampling subspace. Indeed, the coefficient operator
associated to $\pi(\Gamma) \eta$ factors into the continuous wavelet
transform associated to $\eta$, followed by the restriction map
$R_\Gamma$. Conversely,
discretisations give rise to sampling subspaces, as the following proposition
shows.
\begin{prop}
\label{prop:sampl_dwt}
 Let $(\pi,{\cal H}_{\pi})$ be a square-integrable representation of
 $G$, and let $\eta$ be an admissible vector. Assume in addition that
 $\pi(\Gamma) \eta$ is a tight frame of ${\cal H}_{\pi}$ with frame
 constant $c_{\eta}$.
 Then ${\cal H} = V_{\eta}({\cal H}_{\pi})$ is a sampling space, and
 $S = \frac{1}{c_{\eta}} V_{\eta} \eta$ is the associated
 $\Gamma$-sinc-type function for ${\cal H}$.
\end{prop}

\begin{prf}
 Clearly $V_{\eta}
 ({\cal H}_{\pi})$ consists of continuous functions. 
 Using the isometry property of $V_{\eta}$ together with the tight frame
 property of $\pi(\Gamma) \eta$, we obtain for all $f = V_{\eta} \phi
 \in {\cal H}$
 \begin{eqnarray*}
  f & = & V_{\eta} \phi = V_{\eta} \left( \frac{1}{c_\eta}
 \sum_{\gamma \in \Gamma} \langle
 \phi, \pi(\gamma) \eta \rangle \pi(\gamma) \eta \right) \\
 & = & \sum_{\gamma \in \Gamma} \frac{1}{c_\eta}
 V_{\eta} \phi(\gamma) V_{\eta} (\pi(\gamma)
 \eta) = \sum_{\gamma \in \Gamma} f(\gamma) S(\gamma^{-1} \cdot) 
 ~~,
 \end{eqnarray*}
 with convergence in $\| \cdot \|_2$. In order to obtain uniform
 convergence, we find that for $g = V_{\eta} \psi \in {\cal H}$,
 \[ | f(x) - g(x)| = |\langle \phi - \psi, \pi(x) \eta \rangle |
 \le \| \phi - \psi \| ~ \| \eta \| = \frac{1}{c_\eta}
 \| f-g \|_2 \| \eta \| ~~,\]
 i.e., ${\rm L}^2$-convergence entails uniform convergence.
\end{prf}

\begin{rems}
{\bf \underline{1.}}
The original sampling theorem over the reals can be seen to fit into this
setting. If we pick $\eta$ to be the sinc-function, we find that
$V_{\eta} : {\cal H} \to {\rm L}^2(\RR)$ is just the inclusion map,
hence $\eta$ is admissible. Moreover, it is immediately
checked on the Fourier transform side, that $(\lambda_{\RR}(n) \eta)_{n
\in \ZZ}$ is an ONB of ${\cal H}$. \\
{\bf \underline{2.}}
The proposition shows that various results on the relation between discrete
wavelet or Weyl-Heisenberg systems and continuous ones give rise to 
sampling theorems:
For the wavelet case, the underlying group is the $ax+b$-group. A result
by Daubechies \cite{Da} ensures that every wavelet giving rise to a frame 
is in fact an admissible vector, hence we are precisely in the setting of the
proposition. Similarly for discrete Weyl-Heisenberg system, where the 
underlying
group is the Weyl-Heisenberg group (a quotient of the Heisenberg
group by a discrete central subgroup). Here admissibility of the
window function is trivial, as is always the case for irreducible 
square-integrable representations of unimodular groups. Again the 
expansion coefficients
are sampled values of the windowed Fourier transform, which is
the underlying (generalised) continuous wavelet transform.
\end{rems}

Our next aim is to show that, at least for unimodular groups, the
discretisation problem and construction of sampling spaces are
equivalent problems. We first need an auxiliary result ensuring
the existence of sufficiently many convolution idempotents. It
immediately follows from \cite[Theorem 2.3]{Ri}.
\begin{prop} \label{prop:enough_idpt}
 Let $G$ be unimodular, let $\{ 0 \} \not= {\cal H} \subset {\rm L}^2(G)$ be
 closed and leftinvariant. Then ${\cal H}$ contains a nontrivial
 selfadjoint convolution idempotent.
\end{prop}

The following theorem serves various purposes: First of all it
shows that, at least for a large class of settings, the definition of
a sampling space is redundant: Property $(ii)$ follows from property $(i)$.
Secondly it shows that every sampling space can be obtained from the
construction in Proposition \ref{prop:sampl_dwt}, hence the construction
of sampling subspaces and the discretisation problem are (in a somewhat
abstract sense) equivalent. 
\begin{thm} \label{thm:sp_equiv_dwt}
 Assume that $G$ is unimodular.
 Let ${\cal H} \subset {\rm L}^2(G)$ be a leftinvariant closed subspace 
 consisting of continuous functions, and assume that it has property $(i)$
 of a sampling space. Then ${\cal H}$ is a sampling subspace.
 More precisely, there exists a selfadjoint convolution idempotent $S$,
 sucht that $\frac{1}{c_{\cal H}} S$ is the associated sinc-type function,
 and in addition ${\cal H} = {\rm L}^2(G) \ast S$. In particular, 
 \[  
 \forall f \in {\cal H}~,~ \forall \gamma \in \Gamma~~:~~ 
 f(\gamma) = \langle f, \lambda_G(\gamma) S \rangle~~,
 \]
 and thus $\lambda_G(\Gamma) S$ is a tight frame for ${\cal H}$.
 The restriction map $R_{\Gamma}$ is onto iff $\lambda_G(\Gamma) S$
 is an orthonormal basis of ${\cal H}$.
\end{thm}
\begin{prf}
 Let $R_{\Gamma}: {\cal H} \to \ell^2(\Gamma)$ denote the restriction map,
 and define
 \[ S_{\gamma} =  R_{\Gamma}^*(\delta_{\gamma}) ~~,\]
 where $\delta_{\gamma} \in \ell^2(\Gamma)$ is the Kronecker delta
 at $\gamma$. Then $\frac{1}{c_{\eta}} R_{\Gamma}^* R_{\Gamma} = {\rm Id}_{\cal H}$
 shows that
 \begin{equation} \label{eq:sp_exp_1}
 f = \sum_{\gamma \in \Gamma} f(\gamma) \frac{1}{c_{\cal H}}
 S_{\gamma} ~~,\end{equation}
 with convergence in the norm. 
% By polarisation, property $(i)$ gives
% \[ \langle R_{\Gamma}f, R_{\Gamma}g \rangle_{\ell^2} = c_{\cal H} \langle
% f, g \rangle_{{\rm L}^2} ~~.\]
 The orthogonal projection $P: \ell^2(\Gamma) \to R_{\Gamma}({\cal H})$
 is given by $P = \frac{1}{c_{\cal H}} R_{\Gamma} R_{\Gamma}^*$.
 Hence, using polarisation we compute 
 \begin{equation} \label{eqn:samp_sc_1}
  f(\gamma) = \langle R_{\Gamma}f, \delta_{\gamma} \rangle = \frac{1}{c_{\cal H}}
 \langle R_{\Gamma}f, R_{\Gamma} R_{\Gamma}^* \delta_\gamma \rangle = \langle f, R_{\Gamma}^* \delta_{\gamma}
 \rangle = \langle f, c_{\cal H} S_{\gamma} \rangle~~.
 \end{equation}
  
 Next pick a maximal family $({\cal H}_i)_{i \in I}$ of nontrivial pairwise
 orthogonal closed subspaces of the form ${\cal H}_i = {\rm L}^2(G) \ast S_i$,
 where the $S_i$ are
 selfadjoint convolution idempotents in ${\rm L}^2(G)$.
 Then Proposition \ref{prop:enough_idpt} implies that
 \[ {\cal H} = \bigoplus_{i \in I} {\cal H}_i ~~.\]
 Since right convolution
 with $S_i$ is the orthogonal projection onto ${\cal H}_i$,
 equation (\ref{eqn:samp_sc_1}) implies, for all $f \in
 {\cal H}_i$,
 \[ \langle f, S_{\gamma} \ast S_i \rangle = \langle f \ast S_i,
  S_{\gamma} \rangle = \langle f,
  S_{\gamma} \rangle = f(\gamma) = \langle f, \lambda_G(\gamma)
  S_i \rangle ~~.\]
 As a consequence, $ S_{\gamma} \ast S_i = \lambda_G(\gamma) S_i$.
 Letting
 \[ S = \sum_{i \in I} S_i
 \]
 we find that
 \begin{equation} \label{eqn:S_transl}
  S_{\gamma} = \sum_{i \in I} S_{\gamma} \ast S_i 
 = \sum_{i \in I}  \lambda_G(\gamma) S_i = \lambda(\gamma) S ~~.\end{equation}
 (Note that equation (\ref{eqn:S_transl}) implies that the sum defining $S$
 converges in the norm.) The self-adjointness of the $S_i$
 implies the same property for $S$. Moreover, for all $f \in {\cal H}$,
 \[ (f \ast S^*)(x) = \langle f, \lambda_G(x) S \rangle = \langle f, 
 \sum_{i \in I} \lambda_G(x) S_i \rangle = \left( \sum_{i \in I} f \ast S_i
 \right) (x) = f(x) ~~.\]
 Hence ${\cal H}= {\rm L}^2(G) \ast S$. Now (\ref{eq:sp_exp_1}) and 
 (\ref{eqn:S_transl}) shows that for $\frac{1}{c_{\cal H}} S$ to be the
 associated sinc-type function, only the uniform convergence of the
 sampling expansion remains to be shown. For the latter 
 we note once again that
 \[ | f(x) - g(x) | = | \langle f-g, \lambda_G(x) S \rangle |
 \le \| S \|_2 \| f - g \|_2 ~~,\]
 i.e., norm-convergence entails uniform convergence.
 The remaining statements of the theorem are obvious.
\end{prf}

In our treatment of the Heisenberg group we will focus on {\em regular
sampling}, meaning that the sampling set $\Gamma$ is a subgroup.
This setting allows further observations.

\begin{prop}
\label{prop:lattice}
 Let $G$ be a unimodular group, $\Gamma < G$ a discrete subgroup
 and ${\cal H} \subset {\rm L}^2(G)$ a sampling subspace for $\Gamma$.
 Then $\Gamma$ is a lattice, with ${\rm covol}(\Gamma) =
 \frac{1}{c_{\cal H}}$.
\end{prop}

\begin{prf}
 A lattice is by definition a discrete subgroup with finite
 covolume, the latter being defined as ${\rm covol}(\Gamma)
 = \mu_G(A)$, where $A$ is a measurable set of representatives
 of the right coset space $G/\Gamma$. It is straightforward
 to check that the definition does not depend on the choice of $A$. 
 If $f \in {\cal H}$ is any nonzero vector, we can therefore compute
 \[ \| f \|^2 = \int_{A} \sum_{\gamma \in \Gamma} | f(x \gamma) |^2
 d\mu_G(x) = \int_{A} 
 c_{\cal H} \| \lambda_G(x^{-1}) f \|^2 d\mu_G(x)
 = \| f \|^2 c_{\cal H} ~ {\rm covol}(\Gamma) ~~.
 \]
\end{prf}

%For the real case the connection between existence of tight
%frames and the sampling space property is even more obvious:
%\begin{prop}
% Let ${cal H} \subset {\rm L}^2(\RR)$ be a closed translationinvariant
% subspace, and let $\Gamma = r \ZZ$ be a lattice in $\RR$. Then the 
% following are equivalent:
%\begin{enumerate}
% \item[(a)] ${\cal H}$ is a sampling subspace with respect to $\Gamma$.
% \item[(b)] ${\cal H}$ has a tight frame of the form $\lambda_\RR(\Gamma) S$.
% \end{enumerate}
%\end{prop}

\section{The Heisenberg group}

In this section we collect the relevant results concerning the
Heisenberg group.
The Heisenberg Lie algebra ${\mathfrak h}$ is a threedimensional
Lie algebra, with basis $P,Q,Z$ and
defining relations
\[ [P,Q]=Z ~~,\]
all other commutators vanish. It is a one-step nilpotent Lie algebra,
and we equip it with the Campbell-Baker-Hausdorff product to make
it a Lie group $\HH$. It is a unimodular group, with the usual Lebesgue
measure on $\HH \equiv \RR^3$ as Haar measure.
In the coordinates with respect to the above basis,
the product reads
\[ (p,q,t)\ast(p',q',t') = (p+p',q+q',t+t'+(p q' - q p')/2) ~~.\]
The center of $\HH$ is given by $Z(\HH) = \{ (0,0,t) : t \in \RR \}$.
We denote the group of topological automorphisms of ${\HH}$ by
${\rm Aut}(\HH)$.
There exists a family of irreducible, pairwise inequivalent representations 
$\rho_h$ ($h \not=0$)
on ${\rm L}^2(\RR)$, namely the {\bf Schr\"odinger representations}
acting via
\[
\left[ \rho_h(p,q,t) f \right] (x) = e^{2 \pi i h t} e^{ 2\pi i q x}
e^{\pi i h p q} f(x+hp) ~~.
\]
The Schr\"odinger representations do not exhaust the dual of $\HH$,
which in addition contains the characters of the abelian factor group 
$\HH / Z(\HH)$. However, for the decomposition of the regular 
representation of $\HH$, we may concentrate on the Schr\"odinger 
representations. This is a consequence of the Plancherel Theorem, which
we state next. The proof may be found in \cite{Fo}. Some notation
and definitions in connection with direct integrals and Hilbert-Schmidt
operators, which we will use without further comment throughout the
paper, can be found in the appendix.
\begin{thm}
 Define, for $f \in {\rm L}^2(\HH) \cap {\rm L}^1(\HH)$ 
 the {\bf Fourier transform} as the operator field
 \[ {\cal P}(f) := \left( \rho_h(f) \right)_{h \in \RR'} :=
 \left( \int_{\HH}^{wo} f(x) \pi(x) d\mu_{\HH}(x) \right)_{h
 \in \RR'} ~~.\]
 Then ${\cal P}(f)(h) \in {\cal B}_2({\rm L}^2(\RR))$, the space of
 Hilbert-Schmidt operators on ${\rm L}^2(\RR)$, and we have
 the Parseval formula
 \[ \| f \|_{{\rm L}^2}^2 = \int_{\RR'} \|\rho_h(f)\|_{{\cal B}_2}^2
 ~|h| dh .\]
 By continuity this mapping extends to a unitary equivalence
 \[ {\cal P}: {\rm L}^2(\HH) \to \int_{\RR'}^{\oplus} {\cal B}_2({\rm L}^2
 (\RR)) ~|h| dh ~~,\]
 which is the Plancherel transform on $\HH$. It intertwines the
 two-sided regular representation with the direct integral representation
 $\int_{\RR'} \rho_h \otimes \overline{\rho_h} ~|h| dh$
\end{thm}

In the following we shall use $\widehat{f}$ to denote
the Plancherel transform of an ${\rm L}^2$-function $f$. Besides
the decomposition into irreducibles, the Plancherel decomposition 
provides a simple characterization of leftinvariant operators on 
${\rm L}^2(\HH)$.
We remark that the following proposition can be formulated
for general unimodular groups with type I regular
representation \cite{Di}.
\begin{prop} \label{prop:lip_pl}
 Let $T : {\rm L}^2(\HH) \to {\rm L}^2(\HH)$ be a bounded 
 leftinvariant operator. Then there exists a measurable field
 $\widehat{T}_h$ of uniformly bounded operators on ${\rm L}^2(\RR)$ such that
 for all $f \in {\rm L}^2(\HH)$ and almost all $h \in \RR'$,
\[
 \widehat{(T f)}(h) = \widehat{f}(h) \circ \widehat{T}_h~~.
\]
 In particular, given any closed leftinvariant subspace
 ${\cal H} \subset {\rm L}^2(\HH)$,
 the orthogonal projection $P$ onto ${\cal H}$ decomposes
 on the Plancherel transform side into a measurable field 
 $(\widehat{P}_h)_{h \in \RR'}$ of orthogonal projections
 on ${\rm L}^2(\RR)$.
\end{prop}

The proposition motivates the following definitions:
\begin{defn}
For a leftinvariant subspace ${\cal H}$ let $(\widehat{P}_h)_{h \in \RR'}$
denote the associated field of projection operators, and define 
\[ \Sigma({\cal H}) = \{ h \in \RR' :\widehat{P}_h \not= 0 \} \]
and similarly, 
$ \Sigma(g) = \{ h \in \RR' : \widehat{g}(h) \not= 0 \}$,
for $g \in {\rm L}^2(\HH)$. Both are defined only up to a set of measure zero.
We call ${\cal H}$ (resp. $g$) {\bf bandlimited} if 
$\Sigma({\cal H})$ (resp. $\Sigma(g)$) is a bounded set in $\RR$.
In addition, define 
 \[ m(h) :=  {\rm rank}(\widehat{P}_h)~~,\]
 the {\bf multiplicity function} of ${\cal H}$. ${\cal H}$ is called {\bf 
 multiplicity-free} if $m(h) \in \{ 0, 1 \}$, for almost every $h$.
\end{defn}

In the abelian case, say over the reals, $\Sigma({\cal H})$ has a natural 
counterpart, and it plays a central role for the characterisation of sampling
spaces. The multiplicity
function provides more detailed information, which in the abelian
case is superfluous, since $\lambda_\RR$ is multiplicity-free. 

The following theorem characterises admissible vectors for leftinvariant
subspaces via their Plancherel transform. See \cite{Fu} for the proof of 
the general version for unimodular groups with type I regular representation. 
The (rather obvious) part $(iii)$ is not proved there; it provides a 
description of selfadjoint convolution idempotents on the Plancherel 
transform side.

\begin{thm}
\label{thm:adm_vec}
 Let ${\cal H} \subset {\rm L}^2(\HH)$ be leftinvariant, and let
 $(\widehat{P}_h)_{h \in \RR'}$ be the associated field of projection
 operators. 
 \begin{enumerate}
 \item[(i)] $ f \in {\cal H}$ is admissible for ${\cal H}$
 iff $\widehat{f}(h)^*$, restricted to
 ${\widehat{P}_h}({\rm L}^2(\RR))$, is an isometry.
 \item[(ii)] There exist admissible vectors for ${\cal H}$ iff 
 $\int_{\RR'} m(h) |h| dh < \infty$.
 \item[(iii)] The necessary and sufficient condition of $(ii)$ is
 equivalent to the property
 that the operator field $(\widehat{P}_h)_{h \in \RR'}$ is the
 Plancherel transform of a selfadjoint convolution idempotent $S 
 \in {\rm L}^2(\HH)$. Accordingly, ${\cal H} = {\rm L}^2(G) \ast S$.
 \end{enumerate}  
\end{thm} 

To close our survey of the Heisenberg group, we cite a result classifying
the lattices of $\HH$. We associate to such a lattice $\Gamma$
two numbers $d(\Gamma) \in \NN'$, $r(\Gamma) \in \RR^+$ which contain
sufficient information for our purposes. Both parameters can be interpreted
as a measure of the density of $\Gamma$ in $\HH$. We first single out 
a particular family of lattices, which turns out to be exhaustive (up to 
automorphisms of $\HH$).

\begin{defn}
 For any  positive integer $d$ let $\Gamma_d$ be the subgroup
 generated by $P, dQ, Z$. $\Gamma_d$ is a lattice, with 
 \begin{eqnarray*}
 \Gamma_d & = &  (\ZZ \cdot Z) \ast (\ZZ \cdot P) \ast (d \ZZ \cdot Q) \\ 
 & = &  \left\{ (m,dk,\ell+\frac12 d m k) :m,k, \ell \in \ZZ \right\} ~~.
 \end{eqnarray*}
 It is convenient to introduce the {\bf reduced lattice} $\Gamma_d^r$
 which is the subset
 \begin{eqnarray*}
  \Gamma_d^r & = &  (\ZZ \cdot P) \ast (d \ZZ \cdot Q) \\ 
 & = &  \left\{ (m,dk, d m k/2) :m,k\in \ZZ \right\} ~~.
 \end{eqnarray*}
 Note that $\Gamma_d^r$ is not a lattice, not even a subgroup. 
\end{defn}

Let us next give a classification of lattices.
It has been attributed (in more generality) to Maltsev.
Since we were not able to locate a source, we sketch a
short proof for the sake of completeness.
\begin{thm}
 Let $\Gamma$ be a lattice of $\HH$.
 Then there exists a strictly positive integer $d$ and
 $\alpha \in {\rm Aut}(\HH)$ with $\alpha(\Gamma_d) = \Gamma$.
 The integer $d$ is uniquely determined by these properties.
\end{thm}

\begin{prf}
 By \cite[Theorem 5.1.6]{CoGr}, there exist a basis $\widetilde{P},
 \widetilde{Q}, \widetilde{Z}$ of $\mathfrak{h}$ with $\widetilde{Z} \in 
 Z(\HH)$, and $\Gamma = \ZZ \widetilde{Z} \ast \ZZ \widetilde{P}
 \ast \ZZ \widetilde{Q}$. Now $[\widetilde{P},\widetilde{Q}] = 
  \widetilde{P} \widetilde{Q} \widetilde{P}^{-1}\widetilde{Q}^{-1} \in \Gamma
 \cap Z(\HH) = \ZZ \widetilde{Z}$ implies $[\widetilde{P},\widetilde{Q}] = d 
 \widetilde{Z}$ for some $d \in \ZZ$, w.l.o.g. $d\ge 0$ (otherwise
 exchange $\widetilde{Q},\widetilde{P}$). In fact, $d>0$  since
 $\mathfrak{h}$ is not abelian.
 It is immediately checked that the linear isomorphism defined by
 \[ P \mapsto \widetilde{P} ~~,~~dQ \mapsto \widetilde{Q}~~,~~
 Z \mapsto \widetilde{Z} \]
 is in  ${\rm Aut}(\HH)$. That $d$ is unique follows from the fact
 that each automorphism $\alpha$ mapping $\Gamma_d$ to $\Gamma_{d'}$  
 maps $Z$ onto $\pm Z$. From Proposition \ref{prop:aut_2} (a) below
 follows that $\alpha$ leaves the Haar measure of $\HH$ invariant, and this
 implies that ${\rm covol}(\Gamma_d) = {\rm covol}(\Gamma_{d'})$. 
 On the other hand, ${\rm covol}(\Gamma_d) = d$, hence $d=d'$.
\end{prf}

We denote by $d(\Gamma)$ the unique integer $d$ from the theorem. For the 
definition of $r(\Gamma)$ we take $r \in \RR^+$ with
$\Gamma \cap Z(\HH) = r(\Gamma) \ZZ Z$.
 
\section{Main results}
\label{sect:main_results}

Now we can state the main results of this paper. In this section,
${\cal H}$ always denotes a closed, leftinvariant subspace of 
${\rm L}^2(\HH)$, and $\Gamma < \HH$ a lattice. Recall from Theorem
\ref{thm:sp_equiv_dwt} that we may assume ${\cal H} = {\rm L}^2(G) \ast S$,
where $S$ is a selfadjoint convolution idempotent,
and that ${\cal H}$ is a sampling
space iff $\lambda_\HH(\Gamma) S$ is a tight frame of ${\cal H}$.
The main theorem characterises the subspaces admitting tight
frames. 

\begin{thm} \label{thm:main_1}
 \begin{enumerate} 
 \item[(i)] There exists a tight frame of the form $\lambda_{\HH}(\Gamma)\Phi$
 with suitable $\Phi \in {\cal H}$ iff the multiplicity function $m$ associated
 to ${\cal H}$ satisfies
 \begin{equation} \label{eqn:mult_vers_dens}
  m(h) \cdot |h| + m\left( h-\frac{1}{r(\Gamma)} \right)
 \cdot \left|h-\frac{1}{r(\Gamma)}\right| \le \frac{1}{d(\Gamma) r(\Gamma)}~~
 \mbox{(almost everywhere)}.
 \end{equation}
 In particular, $\Sigma({\cal H}) \subset \left[-
 \frac1{d(\Gamma)r(\Gamma)}, \frac1{d(\Gamma)r(\Gamma)} \right]$
 (up to a set of measure zero).
 \item[(ii)] If $\lambda_{\HH}(\Gamma) \Phi$ is a tight frame of ${\cal H}$,
 then $\frac{1}{\sqrt{d(\Gamma)}r(\Gamma)} \Phi$ is admissible for ${\cal H}$.
 \item[(iii)] There does not exist an orthonormal basis of the
 form $\lambda_{\HH}(\Gamma) \Phi$ for ${\cal H}$.
\end{enumerate}
\end{thm}

\begin{rem}
 Note that for $d(\Gamma) > 1$, inequality (\ref{eqn:mult_vers_dens})
 simplifies to
 \begin{equation} \label{eqn:mult_vers_dens_d>1}
  m(h) \cdot |h| \le \frac{1}{d(\Gamma) r(\Gamma)}~~
 \mbox{(almost everywhere)}.
 \end{equation}
\end{rem}

\begin{cor} \label{cor:main_2}
 Assume that the multiplicity function $m$ associated to ${\cal H}$
 is bounded. There exists a lattice $\Gamma$ and a $\Phi \in {\cal H}$
 such that $\lambda_\HH(\Gamma) \Phi$ is a tight frame of ${\cal H}$
 iff ${\cal H}$ is bandlimited.
\end{cor}

The following is a rephrasing for square-integrable representations:
\begin{cor} \label{cor:main_3}
 Let $(\pi, {\cal H}_{\pi})$ be a square-integrable
 representation of $\HH$ with bounded multiplicity. Associate the set
 $\Sigma(\pi) \subset \RR'$
 by picking an admissible vector $\eta \in {\cal H}_{\pi}$
 and letting $\Sigma(\pi) := \Sigma(V_{\eta} \eta)$; up to a null set,
 this is independent of the choice of $\eta$. 
 There exists a lattice $\Gamma$ and a vector $\eta$ such that
 $\pi(\Gamma) \eta$ is a tight frame iff $\Sigma(\pi)$ is bounded.
 Any such $\eta$ is admissible. 
\end{cor}

That bounded multiplicity cannot be dispensed with in the
last corollary is shown by the next result:

\begin{cor} \label{cor:main_4}
 There exists a bandlimited leftinvariant subspace ${\cal H} = 
 {\rm L}^2(G) \ast S$, with a selfadjoint convolution idempotent
 $S \in {\rm L}^2(\HH)$, admitting no tight frame of the form
 $\lambda_\HH(\Gamma) \Phi$.
\end{cor}

With regard to the existence of sampling subspaces, we have:
\begin{cor} \label{cor:main_5} Not every space admitting a tight
 frame of the form $\lambda_{\HH} (\Gamma) S$ is a sampling subspace 
 for $\Gamma$. However, for such a space ${\cal H}$ there exists $\Phi \in
 {\cal H}$ such that $f \mapsto f \ast \Phi^*$ is an isometry on ${\cal H}$,
 mapping ${\cal H}$ onto a sampling space. There does not exist a sampling
 space ${\cal H}$ for which the restriction map $R_\Gamma$ is onto.
\end{cor}

%The following corollary should be contrasted to the fact that
%${\rm L}^2(\RR) =  \bigoplus_{i \in I} {\cal H}_i$, where the ${\cal H}_i$
%are sampling spaces.

%\begin{cor} \label{cor:main_6}
% For a lattice $\Gamma$, define
% \[ {\cal H}_{\Gamma} = \left\{ f \in {\rm L}^2(\HH) : \Sigma(f) \subset
%  \left[-
% \frac1{d(\Gamma)r(\Gamma)}, \frac1{d(\Gamma)r(\Gamma)} \right] \right\}~~.
% \]
% By Theoremm \ref{thm:main_1}, ${\cal H} \subset {\cal H}_{\Gamma}$,
% for every sampling subspace ${\cal H}$ with respect to $\Gamma$.
% Conversely, there exists a family $({\cal H}_i)_{i \in I}$ of pairwise
% orthogonal sampling spaces with ${\cal H}_{\Gamma} =
% \bigoplus_{i \in I} {\cal H}_i$.
%\end{cor}

%The following theorem illustrates that band-limitedness is
%not sufficient for the existence of sinc-type functions.
%The main statement is that there exists a band-limited
%subspace with the property that its element are not characterized
%by their restriction to any lattice, however dense it might
%be. 
%\begin{thm}
% \label{thm:main_5}
% There exists a multiplicity-free closed leftinvariant space
% ${\cal H} \subset {\rm L}^2(\HH)$
% with the following properties: ${\cal H}$ admits $\Gamma$-NTF's, for
% suitably chosen lattices $\Gamma$. In particular ${\cal H}$ is 
% band-limited. On the other hand, for every lattice $\Gamma$,
% the restriction mapping
% \[ {\cal H} \ni f \mapsto (f(\gamma))_{\gamma \in \Gamma}
% \in \ell^2(\Gamma) \]
% is well-defined, bounded, but not injective. 
%\end{thm}

The proofs for these results will be given in Section \ref{sect:proof_main}
below. The following remarks discuss similarities  and differences
to the case of the reals:

\begin{rems} 
\begin{enumerate}
\item The main similarity lies in the notion of bandwidth, and the
 fact that it can be interpreted as inversely proportional to the
 density of the lattice. Note that over $\HH$ the bandwidth restriction
 is much more rigid: The set $\Sigma({\cal H})$
 is contained in a fixed interval, whereas the analog of that set in
 the real case
 can be shifted arbitrarily and still give a sampling subspace.
\item Corollaries \ref{cor:main_4} and \ref{cor:main_5} mark
 important differences between the sampling theories of $\HH$ and 
 $\RR$. None of the counterexamples given in the corollaries 
 has an analog in the real setting. In particular, the question whether a given
 space is a sampling space is much more subtle than deciding whether
 it has a frame. For the first problem, a close inspection of the projection
 operator field $(\widehat{P}_h)_{h \in \RR'}$ is necessary, for the
 second, only the ranks of these operators are needed.
 By contrast, for the reals it can be shown that every subspace
 of ${\rm L}^2(\RR)$ admitting a frame obtained from the action
 of a lattice is already a sampling subspace.
\item While Theorem \ref{thm:main_1} shows that the Plancherel transform
 can be used to characterise sampling spaces and frames, it is not
 clear how it can be generalised to a larger class of locally compact groups.
 Indeed, as far as we are aware, among the
 entities entering the central relation (\ref{eqn:mult_vers_dens}),
 only the multiplicity function $m$ has an abstract interpretation.
 A possible starting point for more general considerations could be
 to study multiplicity-free subspaces and try to come up with criteria
 involving the  natural topology on the dual.
\item We have used lattices as sampling sets simply because
 they are easily accessible. In particular, we have not at all
 exploited the representation theory of the lattices. An alternative
 approach along these lines could provide valuable additional information. 
\end{enumerate}
\end{rems}

\section{Reduction to Weyl-Heisenberg systems}

In this section we start the discussion of normalised tight frames
for leftinvariant subspaces.
On the Plancherel transform side, the space ${\cal H}$ under
consideration decomposes into a direct integral. In this section, we reduce
the complexity of the problem in two ways: We get rid of the
direct integral on the one hand, and the central variable of
the lattice on the other, and are faced with the problem of 
constructing certain normalised tight frames in the fibres, 
arising from the action of the reduced lattice. The latter
problem is equivalent to the construction of Weyl-Heisenberg (super-)frames,
which allows to finish our proof.

\begin{prop}
\label{prop:boil_to_red}
 Let $\Phi \in {\cal H}$ be such that $\lambda(\Gamma) \Phi$ is a normalised
 tight frame of ${\cal H}$.
 Then, for almost every $h \in \Sigma({\cal H})$, the reduced lattice
 satisfies the following condition:
 \begin{equation} \label{eq:fibre_ntf}
  \left( |h|^{1/2}  \rho_h(\gamma) \widehat{\Phi}(h)
 \right)_{\gamma \in
 \Gamma^r} \mbox{ is a normalised tight
 frame of }  {\cal B}_2({\rm L}^2 (\RR)) \circ \widehat{P}_h ~.
 \end{equation}
 Conversely, if both (\ref{eq:fibre_ntf}) (for almost every $h$) and 
 the support condition
 \begin{equation} \label{eq:supp_cond}
  \forall m \in \ZZ \setminus \{ 0 \}~~:~~ \Sigma ({\cal H})
 \cap m + \Sigma({\cal H}) \mbox{ has measure zero}
 \end{equation}
 hold, then $\lambda_G(\Gamma) \Phi$ is a normalised tight frame
 of ${\cal H}$.
\end{prop}

\begin{prf}
 We calculate
 \begin{eqnarray*}
 \sum_{\gamma \in \Gamma} | \langle f, \lambda_\HH
 (\gamma) \Phi \rangle |^2  & = &  \sum_{\gamma \in \Gamma} \left|
 \int_{{\Sigma(f)}} \langle \widehat{f}(h), \rho_h(\gamma) \widehat{\Phi}(h) \rangle
 |h| dh \right|^2 \\
 & = & \sum_{\gamma \in \Gamma^r} \sum_{\ell \in \ZZ}
 \left| \int_{{\Sigma(f)}} e^{-2 \pi i h \ell}  
 \langle \widehat{f}(h), \rho_h(\gamma) \widehat{\Phi}(h) \rangle
 |h| dh \right |^2 \\
 & = & \sum_{\gamma \in \Gamma^r} \int_{{\Sigma(f)}} \left|
  \langle \widehat{f}(h), \rho_h(\gamma) \widehat{\Phi}(h) \rangle
 \right|^2 |h|^2 dh \\
 & = & \int_{{\Sigma(f)}} \sum_{\gamma \in \Gamma^r}  \left|
  \langle \widehat{f}(h), \rho_h(\gamma) |h|^{1/2} \widehat{\Phi}(h) \rangle
 \right|^2 |h| dh ~~. 
\end{eqnarray*}
 Here we used the assumption on $\Sigma(f)$ to apply the Plancherel Theorem
 on $\Sigma(f)$ and thereby discard the summation over $\ell$.
 On the other hand, the tight frame condition together with the Plancherel
 formula for $\HH$ implies that
 \[ \sum_{\gamma \in \Gamma} | \langle f, \lambda(\gamma) s \rangle |^2
 = \int_{{\Sigma(f)}} \|  \widehat{f}(h) \|^2 |h| dh ~~,\]
 and thus
 \begin{equation} \label{eq:int_norms}
 \int_{{\Sigma(f)}} \sum_{\gamma \in \Gamma^r}  \left|
  \langle \widehat{f}(h), \rho_h(\gamma) |h|^{1/2} \widehat{\Phi}(h) \rangle
 \right|^2 |h| dh =  \int_{{\Sigma(f)}} \|  \widehat{f}(h) \|^2 |h| dh~~.
 \end{equation}
 Replacing $f$ by $g$ with $\widehat{g}(h) = \chi_B (h) \widehat{f}(h)$,
 we see that we may replace $\Sigma(f)$ in (\ref{eq:int_norms})
 by any Borel subset $B$. Hence
 the integrands must be equal almost everywhere:
 \begin{equation}
 \label{eq:fibre_norms}
 \sum_{\gamma \in \Gamma^r}  \left|
  \langle \widehat{f}(h), \rho_h(\gamma) |h|^{1/2} \widehat{\Phi}(h) \rangle
 \right|^2 = \| \widehat{f}(h) \|^2   ~~.
 \end{equation}
 Writing an arbitrary $f$ as orthogonal sum of functions $g$ fulfilling
 the initial support condition we see that (\ref{eq:fibre_norms})
 holds for every $f$ and almost every $h \in \RR'$. However, it remains
 to show that the relation holds for all $h$ in a common conull subset, 
 independent of $f$.
 For this purpose we pick a countable dense $\QQ$-subspace ${\cal A}
 \subset {\rm L}^2(\HH)$. Then there exists a conull subset
 $C \subset \RR'$ such that, for all $h \in C$,
 $\{ \widehat{f}(h) : f \in {\cal A} \}$
 is dense in ${\cal B}_2({\rm L}^2(\RR)) \circ \widehat{P}_h$,
 and in addition (\ref{eq:fibre_norms})
 holds for all $f \in {\cal A}$. Now, for every $h \in C$, the
 coefficient map
 \[ \widehat{f}(h) \mapsto \left( \langle \widehat{f}(h), \rho_h(\gamma)
 |h|^{1/2} \widehat{\Phi}(h) \rangle \right)_{\gamma \in \Gamma^r} \]
 is a closed linear operator, by Proposition \ref{prop:basics_frames} (d),
 coinciding with an isometry on a
 dense subset, hence it is an isometry.

 Finally, we note that the argument can be reversed to prove 
 the sufficiency of condition (\ref{eq:fibre_ntf}) under the
 additional assumption (\ref{eq:supp_cond}).
\end{prf}

\section{Weyl-Heisenberg frames}

By a {\bf Weyl-Heisenberg system }${\cal G}(\alpha,\beta,g)$ of ${\rm L}^2(\RR)$ we mean 
a family $(g_{k,m})_{(k,m) \in \ZZ^2}$
\[
 g_{k,m}(x) =  e^{2\pi i \beta k x} g(x+\alpha m) 
\]
resulting from a single function $g \in {\rm L}^2(G)$. A {\bf (normalised,
tight) Weyl-Heisenberg frame} is a Weyl-Heisenberg system which is a (normalised, tight)
frame of ${\rm L}^2(\RR)$. For any $g \in {\rm L}^2(\RR)$,
the operation of the reduced lattice $\Gamma^r$ on $g$ via $\rho_h$ 
gives the system
\[ ( \rho_h(m,dk,  d m k/2) g ) (x) = e^{\pi i h m d k} e^{2 \pi i
  k x} g(x+hm) ~~,\]
hence $\rho_h(\Gamma^r) g$ and the Weyl-Heisenberg system ${\cal G}(h,d,g)$
only differ up to the phase factor $ e^{\pi i h m d k}$. Clearly the
phase factor does not influence any
normalised tight frame or ONB properties of the system, hence we
may and will switch freely between the Weyl-Heisenberg system and the
orbit of the reduced lattice.

The central results concerning Weyl-Heisenberg frames are contained in
the following.
\begin{thm}
\label{lem:Weyl-Heisenberg_density}
 There exists a normalised tight Weyl-Heisenberg frame ${\cal G}(h,d,g)$
 of ${\rm L}^2(\RR)$ iff $|h| d \le 1$.
 For any such frame we have $\| g \|_2^2 = |h| d$.
\end{thm}
\begin{prf}
 The ``only-if''-part is \cite[Corollary 7.5.1]{Gr}. The
 ``if''-part follows from \cite[Theorem 6.4.1]{Gr}, applied to a suitably
 chosen characteristic function.
 The norm equality is due to \cite[Corollary 7.3.2]{Gr}.
\end{prf}

In dealing with subspaces of Hilbert-Schmidt spaces,
we have to consider a more general setting: 
We are interested in normalised tight frames of $\left( 
{\rm L}^2(\RR) \right)^r$
 consisting of vectors 
of the type
\begin{equation} \label{eq:gab_syst_tensor}
  g_{k,m} = ( e^{2\pi i d k x} g^{j}(x+ h_j m) )_{j =1,\ldots,r} = 
 (g^j_{k,m})_{j = 1,\ldots,r}
\end{equation}
where $g = (g^j)_{j = 1,\ldots,r} \in \left( {\rm L}^2(\RR) \right)^r$ 
is suitably chosen,
and $\mathbf{h}=(h_j)_{j} \in \RR^r$ is a vector of nonzero real numbers.
This problem has already been considered by other authors,
see \cite{Ba} and the references therein. Following \cite{Ba},
we call a system of the type (\ref{eq:gab_syst_tensor}) a
{\bf Weyl-Heisenberg superframe}.
The following two lemmata extend the results 
on ${\rm L}^2(\RR)$ to the more general situation. The first one is 
quite obvious and
does not reflect the special structure of Weyl-Heisenberg frames. An alternate
version (for arbitrary frames) is given in \cite{Ba}.

\begin{lemma} \label{lem:NTF_mult_1}
 Let $\mathbf{h} = (h_1,\ldots,h_r)$ and
 $g = (g^j)_{j = 1,\ldots,r} \in  \left( {\rm L}^2(\RR) \right)^r$.
 Then $(g_{k,m})_{k,m \in \ZZ}$, defined as in equation
 (\ref{eq:gab_syst_tensor}), is a normalised tight frame of 
 $\left( {\rm L}^2(\RR) \right)^r$ iff 
 \begin{enumerate}
 \item[(i)] for $j =1,\ldots,r$, ${\cal G}(h_j,d,g^j)$
  is a normalised tight frame of $ {\rm L}^2(\RR)$; and
 \item[(ii)] for $i \not= j$, and for all $f_1,f_2 \in {\rm L}^2(\RR)$,
 \begin{equation} \label{eq:orth_rel_coeff}
 \left\langle \left( \langle f_1, g^j_{m,n} \rangle \right)_{m,n}, 
  \left( \langle f_2, g^i_{m,n} \rangle \right)_{m,n}
 \right\rangle_{\ell^2 (\ZZ \times \ZZ)} = 0 ~~.
 \end{equation}
 i.e., the coefficient operators belonging
 to ${\cal G}(h_j,d,g^j)$ and ${\cal G}(h_i,d,g^i)$ have orthogonal
 ranges in $\ell^2 (\ZZ \times \ZZ)$.
 \end{enumerate}
\end{lemma}
\begin{prf}
 Consider the subspace ${\cal H}_j \subset \left( {\rm L}^2(\RR) \right)^r$
 whose elements are nonzero at most on the $j$th component.
 The necessity of property
 $(i)$ follows immediately from Proposition \ref{prop:basics_frames}(a),
 applied to the ${\cal H}_j$.
 Property $(ii)$ is necessary because the (pairwise orthogonal)
 ${\cal H}_j$ need to have orthogonal
 images in $\ell^2(\ZZ \times \ZZ)$. The converse is clear.
\end{prf}

Necessary and sufficient conditions for the existence of such frames
are given in the next proposition.
\begin{prop} \label{prop:NTF_mult_2}
 Let $(h_j)_{j = 1,\ldots,r}$, $d \in \NN'$ be given.
 \begin{enumerate}
 \item[(a)] There exists a normalised tight frame of $\left(
 {\rm L}^2(\RR) \right)^r$ of the form 
 (\ref{eq:gab_syst_tensor}) iff $d \sum_{j = 1}^{r} |h_j| \le 1$.
 \item[(b)] Assume that $h_j = h$, for all $j = 1,\ldots,r$, and
 $g = (g^j)_{j = 1,\ldots,r}$ is such that (\ref{eq:gab_syst_tensor})
 is a normalised tight frame. Then $g^i \bot g^j$, for $i \not= j$.
\end{enumerate}
\end{prop}

\begin{prf}
 For the necessity in part (a), observe that Lemma \ref{lem:NTF_mult_1}
 together with Theorem \ref{lem:Weyl-Heisenberg_density} yields that
 $\| g^j \|^2 = |h_j|
 d$, and thus $\| g \|^2 = d \sum_{j = 1}^{r} |h_j|$. 
 Now Proposition \ref{prop:basics_frames} (c) entails the desired
 inequality.

 The proof for sufficiency is a slight modification of a construction
 given by Balan \cite[Example 13]{Ba}. Define $c_i = \sum_{j=1}^i |h_j|$,
 and let $g^i = \sqrt{d} \chi_{[c_{i-1},c_i]}$. Given $f = (f^i)
 \in \left( {\rm L}^2(\RR) \right)^r$, we compute
 \begin{eqnarray*}
 \langle f,  g_{k,m} \rangle & = & \sum_{i = 1}^r \langle
 f^i, g_{k,m}^i \rangle \\
 & =  & \sum_{i = 1}^r \sqrt{d}
  \int_{c_{i-1}}^{c_{i}} e^{- 2 \pi i m d x}
 f^i(x+ \beta_i k) dx \\
 & = & \sqrt{d} \int_0^{1/d} e^{-2 \pi i m d x} H_k(x) dx ~~,
 \end{eqnarray*}
 where
 \[ H_k (x) = \left\{ \begin{array}{ll} f^i(x-h_i k) & 
 x \in [c_{i-1},c_i] \\ 0 & \mbox{elsewhere} 
 \end{array} \right. ~~.\]
 Fixing $k$, we compute
 \begin{eqnarray*}
  \sum_{m \in \ZZ} | \langle f, g_{k,m} \rangle |^2 & = & 
 \sum_{m \in \ZZ} d \left|  \int_0^{1/d} e^{-2 \pi i m d x} H_k(x) dx
 \right|^2 \\ & = &  \int_0^{1/d} | H_k(x)|^2 dx \\
 & = & \sum_{i = 1,\ldots,r} \int_{c_{i-1}}^{c_{i}} |f^{i}
 (x-h_i k)|^2 dx~~.
 \end{eqnarray*}
 Since the $h_i \ZZ$-translates of $[c_{i-1},c_i]$ tile
 $\RR$, summing over $k$ yields the desired normequality.
 This closes the proof of (a).

 For the proof of (b), pick $f_1, f_2 \in L^{\infty} (\RR)$
 with supports in $[0,|h|]$.
 Then we calculate 
 \begin{eqnarray*}
 \lefteqn{\sum_{m,k \in \ZZ} \langle f_1, g^j_{m,k} \rangle \overline{
 \langle f_2,
 g^i_{m,k \in \ZZ} \rangle } = }\\
 & = & \sum_{m \in \ZZ} \int_{\RR} \left(
 \sum_{k \in \ZZ} \langle f_1 g^j_{m,0}, e^{2 \pi i d k \cdot} \rangle
 e^{2 \pi i d k x} \right)  \overline{f_2(x)} g^i (x+ h m) dx \\
 & = & \sum_{m \in \ZZ} d^{-1} \int_{0 }^{|h|}
 f_1(x) \overline{f_2(x)} \overline{g^j(x+h m)} 
 g^i(x+h m) dx \\
 & = & d^{-1} \int_0^{|h|}  \left( \sum_{m \in \ZZ}
 \overline{g^j(x+h m)} g^i(x+h m) \right) f_1(x) \overline{f_2(x)} dx 
 \end{eqnarray*}
 Here the Fourier series 
\[
  \sum_{k \in \ZZ} \langle f_1 g^j_{m,0}, e^{2 \pi i d k \cdot} \rangle
 e^{2 \pi d k x} = d^{-1} f_1(x) g^j(x+ h m)
\]
 is valid on $[0, |h| ]$, at least in the ${\rm L}^2$-sense,
 because of $|h| \le d^{-1}$, the latter being a consequence of Theorem
 \ref{lem:Weyl-Heisenberg_density}. Now, for arbitrary $f_1, f_2$, the
 scalar product we started with has to be zero, whence we
 obtain for almost every $x \in [0,|h|]$, 
 \[
 \sum_{m \in \ZZ} g^i(x+ h m)
 \overline{g^j(x+h m)} = 0~~. 
 \]
 Integrating over $[0,|h|]$ and applying Fubini's theorem
 yields $\langle g^i,g^j \rangle = 0$. 
\end{prf}

\begin{rem} \label{rems:constr_frame}
 Note that the vectors $(g^i)_{i = 1,\ldots,r}$ constructed in the
 proof of part $(a)$  
 depend measurably on $\mathbf{h}$, i.e., if we let $(g^i_{\mathbf{h}})$
 be the vector of functions constructed from $\mathbf{h}$,
 then $(x,\mathbf{h}) \mapsto (g^i_{\mathbf{h}}(x))_{i=1,\ldots, r}$ is
 a measurable mapping.
% \underline{2.} Clearly, the $(g^i)$ can be replaced by any
% $\widetilde{g}^i$ with $|\widetilde{g}^i| = |g^i|$.
\end{rem}

\section{Proofs of the main results}
\label{sect:proof_main}

The general proof strategy consists in explicit calculation for the
$\Gamma_d$ and then transferring the results to arbitrary lattices by
the action of ${\rm Aut}(\HH)$.
For this purpose we need a more detailed description of ${\rm Aut}(\HH)$
and its action on the Plancherel transform side. Most of the
results are standard, and we only sketch the proofs.

\begin{prop}
\label{prop:aut_1}
\begin{enumerate}
\item[(a)] For $r>0$ let $\alpha_r(p,q,t) := (\sqrt{r}p,\sqrt{r}q,rt)$.
 Then $\alpha_r \in {\rm Aut}(\HH)$.
 In addition, $\alpha_{inv}: (p,q,t) \mapsto (q,p,-t)$ defines
 an involutory automorphism of $\HH$.
 \item[(b)] Each $\alpha \in {\rm Aut}(\HH)$ can be written uniquely
 as $\alpha = \alpha_r \alpha_{inv}^i \alpha'$,
 where $r \in \RR'$, $i \in \{ 0, 1 \}$ and $\alpha'$ leaves the center
 of $\HH$ pointwise fixed.
 \item[(c)] Suppose that $\alpha(\Gamma_d) = \Gamma$ for some $d$,
 $\alpha$, and let $\alpha = \alpha_s \alpha_{inv}^i \alpha'$ be the
 decomposition from part (b). Then $r(\Gamma) = s$.
\end{enumerate}
\end{prop}
\begin{prf}
 For parts $(a)$, $(b)$ see \cite[Theorem 1.22]{Fo}. Part $(c)$ follows
 directly from the definition of $r(\Gamma)$ and the fact that 
 $\alpha'$ and $\alpha_{inv}$ map every discrete subgroup of $Z(\HH)$
 onto itself. 
\end{prf}

Next let us consider the action on the Fourier transform side.
\begin{prop} \label{prop:aut_2}
\begin{enumerate}
 \item[(a)] Define $\Delta: {\rm Aut}(\HH) \to \RR^+$ by
 \[ \Delta(\alpha) = \frac{\mu_\HH(\alpha(B))}{\mu_\HH(B)} ~~,\]
 where $B$ is a measurable set of positive Haar measure. $\Delta$
 does not depend on the choice of $B$, and it is a continuous group
 homomorphism. For $\alpha = \alpha_r \alpha_{inv}^i \alpha'$ as
 in \ref{prop:aut_1}(b), $\Delta(\alpha) = r^2$.
\item[(b)] For $\alpha \in {\rm Aut}(\HH)$, let ${\cal D}_{\alpha} :
 {\rm L}^2(\HH)
 \to {\rm L}^2(\HH)$ be defined as $({\cal D}_\alpha f)(x) :=
 \Delta(\alpha)^{1/2} f(\alpha(x))$. This defines a unitary operator. 
\item[(c)]
 Let ${\cal H} \subset {\rm L}^2(G)$ be a closed, leftinvariant
 subspace with multiplicity function $m$. Then $\widetilde{{\cal H}} = 
 {\cal D}_{\alpha}({\cal H})$ is closed and leftinvariant as well. Let
 $\widetilde{m}$ denote the multiplicity function related to
 $\widetilde{{\cal H}}$. If $\alpha = \alpha_r
 \alpha_{inv}^i \alpha'$ then $\widetilde{m}$ satisfies
 \begin{equation} \label{eq:auto_mult}
 \widetilde{m}(h) = m( (-1)^i r^{-1} h) ~\mbox{(almost everywhere)}~~.
 \end{equation}
\item[(d)] Let $\Gamma$ be a lattice, $\alpha \in {\rm Aut}(\HH)$ 
 such that $\alpha(\Gamma_d) = \Gamma$. Let ${\cal H} \subset
 {\rm L}^2(\HH)$ be a closed, leftinvariant subspace. Then
 $\lambda_\HH(\Gamma) \Phi$ is a normalised tight frame (an ONB)
 for ${\cal H}$
 iff $\lambda_\HH (\Gamma_d) ({\cal D}_{\alpha}\Phi)$ is a
 normalised tight frame (an ONB) for ${\cal D}_{\alpha}({\cal H})$. 
\end{enumerate}
\end{prop}
\begin{prf}
 Parts (a) and (b) are standard results concerning the action of automorphisms
 on locally compact groups, see \cite{HeRo}. The explicit formula
 for $\Delta(\alpha)$ follows from the fact that every automorphism
 leaving the center invariant factors into an inner
 and a symplectic automorphism \cite[Theorem 1.22]{Fo}; both do not
 affect the Haar measure.

 For part (c), we first note
 that by the Stone-von Neumann theorem \cite[Theorem 1.50]{Fo}, 
 any automorphism $\alpha'$ keeping the center pointwise fixed acts
 trivially on the dual of $\HH$. Hence, 
 \[ \left( {\cal D}_{\alpha'} f \right)^{\wedge}(h) = U_{\alpha',h} \circ
 \widehat{f}(h)  \circ U_{\alpha',h}^* ~~,\]
 where $U_{\alpha',h}$ is a unitary operator on ${\rm L}^2(\RR)$.
 Hence the action of $\alpha'$
 does not affect the multiplicity function, and from now on, we only
 consider $\alpha = \alpha_r \alpha_{inv}^i$. In this case, letting
 \[ (D_r f) (x) = r^{1/2} f(rx) ~~,\]
 we obtain by straightforward computation that
 \begin{equation} \label{eq:dual_action}
  \left( {\cal D}_{\alpha} f \right)^{\wedge}(h) = 
 r^{-1} \cdot D_r \circ \widehat{f}((-1)^i r^{-1} h)
 \circ D_r^* ~~.\end{equation}
 This immediately implies (\ref{eq:auto_mult}).

 To prove (d), observe that the unitarity of ${\cal D}_{\alpha}$ implies
 that ${\cal D}_{\alpha} \left( \lambda_{\HH}(\Gamma) \right)$
 is a normalised tight frame of $ {\cal D}_{\alpha} ({\cal H})$,
 and check the equality
 \[ {\cal D}_{\alpha} (\lambda_\HH(x) S) = \lambda_\HH (\alpha^{-1}(x))
  ({\cal D}_{\alpha} S) ~~.\]
\end{prf}

\noindent
{\bf Proof of Theorem \ref{thm:main_1}.} 
We first prove the theorem for the case $\Gamma = \Gamma_d$. Writing 
 \[ \widehat{\Phi}(h) = \sum_{i \in I_h} \varphi_i^h \otimes \eta_i^h ~~,
 \]
we find by Proposition \ref{prop:boil_to_red}, that for almost every $h$,
$(\rho_h(\gamma) \circ |h|^{1/2} \widehat{\Phi}(h))_{\gamma \in \Gamma}$ has to
 be a normalised tight frame of ${\rm L}^2(\RR) \circ \widehat{P}_h$,
or equivalently, that the vector $(\varphi_i^h)_{i=1,\ldots,m(h)}$
generates a Weyl-Heisenberg superframe of $\left({\rm L}^2(\RR)\right)^{m(h)}$,
for $\mathbf{h}=(h,\ldots,h)$. (For the equivalence, confer 
(\ref{eqn:equiv_tf}) in the appendix.)
 Then \ref{lem:NTF_mult_1}
(a) implies that ${\cal G}(h,d,|h|^{1/2} \varphi_i^h)$ is a normalised
tight frame of ${\rm L}^2(\RR)$. In particular, Theorem 
\ref{lem:Weyl-Heisenberg_density} entails 
\begin{equation} \label{eq:phi_norm}
 \| \varphi_i^h \|^2 = d ~~,\end{equation}
as well as $\Sigma({\cal H}) \subset \left[-
 \frac1{d}, \frac1{d} \right]$.
Moreover, Proposition \ref{prop:NTF_mult_2}(b)
entails that the $\varphi_i^h$ are
mutually orthogonal (for $h$ fixed). This shows that
\[ \widehat{\Phi}(h)^* = \sum_{i \in I_h} \eta_i^h \otimes \varphi_i^h \]
is $\sqrt{d}$ times an isometry on $\widehat{P}_h({\rm L}^2(\RR))$,
and thus $\frac{1}{\sqrt{d}} \Phi$ is admissible for ${\cal H}$,
by Theorem \ref{thm:adm_vec} (ii).
This proves part $(ii)$ of the theorem. If $d > 1$, the support
condition (\ref{eq:supp_cond})
in Proposition \ref{prop:boil_to_red} is fulfilled.
Hence Proposition \ref{prop:NTF_mult_2} (a), applied to 
$\mathbf{h}=(h,\ldots,h)$,
shows that (\ref{eqn:mult_vers_dens_d>1}) is necessary
and sufficient for the existence of a normalised tight frame for ${\cal H}$.
(Note that by Remark \ref{rems:constr_frame},  \ref{prop:NTF_mult_2} (a)
provides a {\em measurable} field of operators.)

The case $d=1$ requires a somewhat more involved argument. Assume that
$\lambda_\HH(\Gamma) \Phi$ is a normalised tight frame, and let $f \in {\cal H}$.
Condition (\ref{eq:fibre_ntf}) from Proposition \ref{prop:boil_to_red}
yields
\begin{equation} \label{eq:alias_1}
 \| f \|^2 = \int_0^1 \left(
 \sum_{\gamma \in \Gamma^r} \left| \langle \widehat{f}(h),
 \rho_h(\gamma) \widehat{\Phi}(h) \rangle \right|^2 |h|^2
 +  \left| \langle \widehat{f}(h-1),
 \rho_{h-1}(\gamma) \widehat{\Phi}(h-1) \rangle \right|^2 |h-1|^2 \right)
 dh ~~.
 \end{equation}
On the other hand,
\begin{eqnarray} \nonumber
\| f \|^2 & = & \sum_{\gamma \in \Gamma^r} \sum_{\ell \in \ZZ}
 \left| \int_{0}^1 e^{-2 \pi i h \ell}  
 \left( \langle \widehat{f}(h), \rho_h(\gamma) \widehat{\Phi}(h) \rangle
 |h| +  \langle \widehat{f}(h-1), \rho_{h-1}(\gamma) \widehat{\Phi}(h-1)
 \rangle |h-1| \right)
 dh \right |^2 \\ \label{eq:alias_2} 
& = &\int_0^1 \left(
 \sum_{\gamma \in \Gamma^r} \left| \langle \widehat{f}(h),
 \rho_h(\gamma) \widehat{\Phi}(h) \rangle |h|
 + \langle \widehat{f}(h-1),
 \rho_{h-1}(\gamma) \widehat{\Phi}(h-1) \rangle |h-1| \right|^2 \right)
 dh ~~.
 \end{eqnarray}
As in the proof of Proposition \ref{prop:boil_to_red}, the fact that
the two equations hold for all $f \in {\cal H}$ allows to equate the
integrands of (\ref{eq:alias_1}) and (\ref{eq:alias_2}). But this implies
the orthogonality of the coefficient families:
\[ \left\langle  \left( \langle \widehat{f}(h),
 \rho_h(\gamma) \widehat{\Phi}(h) \rangle \right)_{\gamma \in \Gamma^r},
 \left( \langle \widehat{f}(h-1), 
 \rho_{h-1}(\gamma) \widehat{\Phi}(h-1) \rangle \right)_{\gamma \in \Gamma^r}
 \right\rangle_{\ell^2(\Gamma^r)} = 0 ~~.\]
Plugging this fact, together with condition (\ref{eq:fibre_ntf}) from
Proposition \ref{prop:boil_to_red}, into Proposition \ref{lem:NTF_mult_1},
we finally realise that the system 
 \[ \left( (\rho_{h-1}(\gamma) |h-1|^{1/2}
\varphi_1^{h-1},..., \rho_{h-1}(\gamma) |h-1|^{1/2} \varphi_{m(h-1)}^{h-1}, 
\rho_h (\gamma) |h|^{1/2} \varphi_1^h,\ldots, \rho_h(\gamma) |h|^{1/2}
\varphi_{m(h)}^h) \right)_{\gamma \in \Gamma^r}
\]
has to be a normalised
tight frame of $\left( {\rm L}^2(\RR) \right)^{m(h)+m(h-1)}$.
An application of Proposition \ref{prop:NTF_mult_2} (a) with
 $\mathbf{h}= (h-1,\ldots,h-1,h,
\ldots,h)$ yields that such a frame exists iff $m(h) |h| +
m(h-1)|h-1| \le 1$. This shows the necessity of (\ref{eqn:mult_vers_dens}).
The sufficiency is obtained by running the proof backward;
the measurability of the constructed operator field is again ensured by Remark 
\ref{rems:constr_frame}.

For the proof of $(iii)$ we need to show, by \ref{prop:basics_frames}(c), that
 $\| \Phi \| < 1$, for 
every $\Phi$ for which $\lambda_{\HH}(\Gamma) \Phi$ is a normalised tight
frame. Recalling that
\[ |h| \| \widehat{\Phi} (h) \|_{{\cal B}_2}^2 = |h| m(h) d~~,\]
and using the fact that the inequality $ m(h) |h| d+  m(h-1) |h-1| d \le 1$
is strict almost everywhere (say, for $h$ irrational) we can
estimate
\begin{eqnarray*}
 \|\Phi \|^2 & = & \int_{-1}^1 \| \widehat{\Phi}(h) \|_{{\cal B}_2}^2 |h|dh \\
 & = & \int_0^1 m(h)|h| d + m(h-1) |h-1| d~ dh \\
 & < & 1~~.
\end{eqnarray*}
This closes the proof for $\Gamma= \Gamma_d$. For 
$\Gamma = \alpha(\Gamma_d)$, write $\alpha = 
\alpha_{r(\Gamma)} \alpha_{inv}^i \alpha'$ as in Proposition
\ref{prop:aut_1} (c). By \ref{prop:aut_2} (d),  we may consider $\Gamma_d$ and
$\widetilde{{\cal H}} = {\cal D}_{\alpha}({\cal H})$
instead of $\Gamma$ and ${\cal H}$. Parts $(iii)$ immediately follows
from this observation. For part $(i)$, we find by Proposition
\ref{prop:aut_2}(c) that
the associated multiplicity function $\widetilde{m}$ fulfills
$\widetilde{m}(h) = m((-1)^i r(\Gamma)^{-1} h)$. Hence, 
(\ref{eqn:mult_vers_dens}) for $\Gamma_d,\widetilde{{\cal H}}$
becomes 
\[
 m((-1)^i r(\Gamma)^{-1} h) |h| + m((-1)^ir(\Gamma)^{-1}(h-1)) |h-1|
 \le \frac{1}{d} ~~ \mbox{(almost everywhere)}
\]
which after dividing both sides by $r(\Gamma)$ and passing to the
variable $\widetilde{h}= (-1)^i r(\Gamma)^{-1}h$ is the desired
inequality (\ref{eqn:mult_vers_dens}).

For part $(ii)$ it remains to show that $\frac{1}{\sqrt{d(\Gamma)}r(\Gamma)}
\widehat{\Phi}(h)$ is an isometry on $\widehat{P}_h({\rm L}^2(\RR))$,
by Theorem \ref{thm:adm_vec}(b). Clearly
the only problem is showing that the normalisation is correct.
Part $(ii)$ for $\Gamma_{d(\Gamma)}$ implies that 
$\frac{1}{\sqrt{d(\Gamma)}}
({\cal D}_{\alpha} \Phi)^{\wedge}(h)$ is a partial isometry, for almost
every $h$. Now relation (\ref{eq:dual_action}) implies that 
$\frac{1}{\sqrt{d(\Gamma)}r(\Gamma)} \widehat{\Phi}(h)$ is a 
partial isometry as well.
\hfill $\Box$

\noindent
{\bf Proof of Corollary \ref{cor:main_2}.} The assumptions
imply that $m(h)|h| \le c$, for all $h \in \RR'$, and $c$ a
constant. Hence picking $s \ge \frac{2c}d$ and
defining $\Gamma = \alpha_s (\Gamma_d)$ ensures that
(\ref{eqn:mult_vers_dens}) is fulfilled.
\hfill $\Box$

\noindent
{\bf Proof of Corollary \ref{cor:main_3}.} For two admissible vectors
$\eta,\eta'$, the spaces $V_{\eta}({\cal H}_{\pi})$ and 
$V_{\eta'}({\cal H}_{\pi})$
carry equivalent subrepresentations of $\lambda_G$. The intertwining
operator $T$ between these subspaces decomposes on the Plancherel transform
side, by Proposition \ref{prop:lip_pl}, and maps $V_{\eta} \eta$
to $V_{\eta'} \eta'$, and thus $\Sigma(V_{\eta} \eta)=\Sigma(V_{\eta'} \eta')$.
This shows the well-definedness of $\Sigma(\pi)$. The remaining statements
are then obtained by transferring the corresponding results from Theorem
\ref{thm:main_1} and Corollary \ref{cor:main_2} back to 
${\cal H}_{\pi}$ via $V_{\eta}^{-1}$. \hfill $\Box$

\noindent
{\bf Proof of Corollary \ref{cor:main_4}.} Pick any measurable
function $m:[-1,1] \to \NN'$ such that $h \mapsto m(h)|h|$ is
integrable but unbounded. There exists a closed, leftinvariant
space ${\cal H}$ with multiplicity function $m$, by Lemma \ref{lem:enough_m}
below. ${\cal H}$ is of the desired form, but violates
(\ref{eqn:mult_vers_dens}), for all lattices $\Gamma$.  \hfill $\Box$

\noindent
{\bf Proof of Corollary \ref{cor:main_5}.} To give an example
proving the first statement, let $\Gamma = \Gamma_d$; using the
appropriate $\alpha \in {\rm Aut}(\HH)$ the argument
can be adapted to suit any other lattice.
For $h \in \left[ 0,
\frac{1}{d} \right]$, define
\[ \eta^h = \frac{1}{\sqrt{h/2}} \chi_{[0,h/2]} ~~.\]
 and $S \in {\rm L}^2(\HH)$ with
$\widehat{S}(h) = \eta_h \otimes \eta_h$. Then $S$ is a selfadjoint
convolution idempotent, and ${\cal H} = {\rm L}^2(\HH) \ast S$ 
has a tight frame of the form $\lambda_\HH(\Gamma) \Phi$. However,
for ${\cal H}$ to be a sampling space, $\lambda_\HH(\Gamma) S$
must be a tight frame, and condition (\ref{eq:fibre_ntf}) implies
that ${\cal G}(h,d,\eta_h)$ is a tight frame of ${\rm L}^2(\RR)$, for
 almost every $h$.
But $\chi_{[h/2,h]}$ has disjoint support with all elements
of that system, hence ${\cal G}(h,d,\eta_h)$ is not even total. 

The second statement is obvious from Theorem \ref{thm:main_1} (ii)
and Proposition \ref{prop:sampl_dwt}. The last statement follows
from Theorem \ref{thm:main_1} (iii). \hfill $\Box$

%\noindent
%{\bf Proof of Corollary \ref{cor:main_6}.} We only prove the corollary
% for $\Gamma=\Gamma_d$. Let $I = \{ \pm 1 \}
% \times \ZZ^2$, and for $I \ni i = {s,k,m}$, define
% \[ \eta_{k,m}^h(x) = \frac{1}{\sqrt{|h|}} \chi_[k|h|,(k+1)|h|](x)
% e^{2 \pi i k x / h} ~~,\]
% and 
% \[ \widehat{S_i}(h) = \chi{\left[0,d(\Gamma)^{-1} r(\Gamma)^{-1} \right]}
% (s h) \eta_{k,m}^h \otimes \eta_{k,m}^h ~~.\]
% Then the $\widehat{S_i}$ are fields of rank.one projection operators
% with integrable trace, hence they are the Plancherel transforms of
% selfadjoint convolution idempotents $S_i$. ${\cal H}_i = 
% {\rm L}^2(\HH) \ast S$ is a sampling subspace: The condition 
% (\ref{eq:fibre_ntf})  from Proposition 
% \ref{prop:boil_to_red} is fulfilled by Remark \ref{rems:constr_frame}
% \underline{2.} Since in addition the support condition
% (\ref{eq:supp_cond}) holds, $\lambda_\HH(\Gamma) S_i$
% is a tight frame. The fact that
% $ \eta_{k,m}^h \bot \eta_{k',m'}^h $ for $(m,k) \not= (m',k')$
% implies that ${\cal H}_i \bot {\cal H}_{i'}$. Finally,
% the totality of the ${\cal H}_i$ in ${\cal H}_{\Gamma}$ is
% easily established. 
% \hfill $\Box$

\section{A concrete example}

In this section we explicitly compute a sinc-type function for
$\Gamma=\Gamma_1$. The construction proceeds backwards, starting
on the Plancherel transform side by giving a field of rank-one projection
operators fulfilling the additional requirements for the sampling
space property. Fourier inversion yields the sinc-type function $S$.
As a consequence, the sampling space is given as ${\rm L}^2(\HH) \ast S$.
In order to minimise tedium, we have
drastically shortened some of the more straightforward calculations.
The three steps carry out the abstract program
developed above.

\noindent
{\bf 1. Construction on the Plancherel transform side:} For $h \in
 [-0.5,0.5]$ let $\eta_h =  |h|^{-1/2} \chi_{[-|h|/2,|h|/2]}$, and
 \[ \widehat{S}(h) = \eta_h \otimes \eta_h ~~,\]
 and let $\widehat{S}$ be zero outside of $[-0.5,0.5]$.
 $\widehat{S}$ is a measurable field of rank-one projection operators, with
 integrable trace, hence has an inverse image $S \in {\rm L}^2(\HH)$  which
 is a selfadjoint convolution idempotent.
 Moreover, it is straightforward to check that 
 $\rho_h(\Gamma_1^r) |h|^{1/2} \eta_h = \rho_h(\Gamma_1^r)
 \chi_{[-0.5,0.5]}$ is a normalised tight frame of ${\rm L}^2(\RR)$,
 (compare the proof of Proposition \ref{prop:NTF_mult_2} (a)).
 Hence, by Proposition \ref{prop:boil_to_red}, $\lambda_\HH(\Gamma) S$
 is a normalised tight frame of ${\cal H} = {\rm L}^2(\HH) \ast S$, and
 ${\cal H}$ is a sampling space.

\noindent
{\bf 2. Plancherel inversion:} We use the inversion formula
\[ f(x) = \int_{\widehat{G}} {\rm trace}(\widehat{f}(\sigma) \sigma(x)^*)
 d\nu_G(\sigma) ~~,\]
 for all $f \in {\rm L}^2(G)$, for which the formula makes sense; i.e.,
 for all $f$ such that $\widehat{f}(\sigma)$ is trace-class ($\nu_G$-almost
 everywhere) and in addition
 $\int_{\widehat{G}} {\rm trace} |\widehat{f}(\sigma)| d\nu_G(\sigma) <
 \infty$.
 The formula was proved by Lipsman \cite{Li}.
 In our concrete
 situation, it is immediately checked that the integrability condition
 is fulfilled, and we obtain
 \begin{eqnarray} 
  S(p,q,t) & = & \int_{-0.5}^{0.5} \langle \eta_h, \rho_h(p,q,t) \eta_h
  \rangle |h| dh \nonumber \\ \label{eq:four_inv_conc}
  & = & \int_{-0.5}^{0.5} e^{-2 \pi i h (t+pq/2)}
  \int_{-\frac{|h|}2}^{\frac{|h|}2} e^{-2\pi i q x} \chi_{[-|h|/2,|h|/2]}
 (x+hp) dx ~ dh ~~.
 \end{eqnarray}

\noindent
{\bf 3. Explicit integration:} Let $\widetilde{S}(p,q,h)$ denote
 the inner integral. In the following, we assume that $q \not= 0$
 and $p \ge 0$. The missing values will be obtained by taking limits
 (for $q=0$) and reflection (for $p<0$). Observe further that $S(p,q,t) = 0$
 for $|p|>1$, hence we will use $|p| \le 1$ wherever we may need it.
 Integration yields 
 \[ \widetilde{S}(p,q,h) = \left\{ \begin{array}{ll}
  \frac{\displaystyle e^{2 \pi i q |h|/2} - e^{-2 \pi i q (|h|/2 - hp)}}{
 \displaystyle 2 \pi i q} & 
 h \ge 0 \\ \frac{\displaystyle
 e^{2 \pi i q (|h|/2 + hp)} - e^{ -2\pi i q |h|/2}}{\displaystyle 2 \pi i q}
 & h < 0 \end{array} \right. ~~,
 \]
 After plugging this into (\ref{eq:four_inv_conc}) and integrating, 
 straightforward  simplifications lead to 
 \[ S(p,q,t) = \frac{1}{2 \pi q} \left( \frac{ \cos \left( \pi(t+(p-1)q/2)
 \right)
 -1}{\pi (t+(p-1)q/2)} - \frac{ \cos \left( \pi( t-(p-1)q/2 ) \right)
 -1}{\pi (t-(p-1)q/2)} \right) ~~.\]
 In order to further simplify this expression, we use the relation
 \[ 
 \frac{\cos (\pi \alpha) - 1}{\pi \alpha} = -\frac{\pi \alpha}2 
 {\rm sinc}^2 \left( \frac{\alpha}2 \right)
 \]
 by which means we finally arrive at 
 \begin{equation} \label{eq:S_expl}
  S(p,q,t) = \frac{1}{4} \left( \left( \frac{t}q + \frac{1-p}2 \right)
 {\rm sinc}^2 \left( \frac{t}2 + \frac{1-p}4 q \right) -
 \left( \frac{t}q - \frac{1-p}2 \right)
 {\rm sinc}^2 \left( \frac{t}2 - \frac{1-p}4 q \right) \right)~~.
 \end{equation}
 For $p<0$ we use that $S(p,q,t) = S^*(p,q,t) = S(-p,-q,-t)$.
 It turns out that
 replacing $p$ by $|p|$ in (\ref{eq:S_expl})
 is the only necessary adjustment for the formula to hold in the general
 case. Finally, sending $q$ to $0$ allows to compute the values
 $S(p,0,t)$, since $S$ is continuous. The following theorem summarises
 our calculations:

 \begin{thm}
 Define $S \in {\rm L}^2(\HH)$ by
 \[ S(p,q,t) = \left\{ \begin{array}{ll} 0 & \mbox{ for } |p|>1 \\
 \frac{1}{4} \left[ \left( \frac{t}q + \frac{1-|p|}2 \right)
 {\rm sinc}^2 \left( \frac{t}2 + \frac{1-|p|}4 q \right) \right.& \\ -
 \left. \left( \frac{t}q - \frac{1-|p|}2 \right)
 {\rm sinc}^2 \left( \frac{t}2 - \frac{1-|p|}4 q \right) \right]& 
 \mbox{ for } |p| \le 1, q \not= 0 \\
  \frac{1-|p|}4 (2 {\rm sinc}(t) - {\rm sinc}^2(t/2)) & \mbox{ for }
 |p| \le 1, q = 0 \end{array} \right.
\]
 Let ${\cal H} \subset {\rm L}^2(\HH)$ be the leftinvariant closed
 subspace generated by $S$. Then ${\cal H}$ is a sampling space for the
 lattice $\Gamma_1$, with $c_{\cal H}=1$, and
 $S$ the associated sinc-type function. $\lambda_{\HH}(\Gamma_1) S$
 is a normalised tight frame, but not an orthonormal basis of ${\cal H}$,
 because of $\| S \|_2 = \frac{1}2$.
 \end{thm}

\begin{appendix}

\section{Hilbert-Schmidt operators and direct integrals}

 In this section we collect a few technical details concerning
Hilbert-Schmidt operators and direct integral Hilbert spaces.
If ${\cal K}$ is a Hilbert space, then
${\cal B}_2({\cal K})$ denotes the space of bounded operators $T$ for
which $T^* T$ is trace-class. ${\cal B}_2({\cal K})$ is a Hilbert
space, with scalar product $\langle S, T \rangle = {\rm tr}(T^*S)$-
The finite-rank operators are dense in ${\cal B}_2({\cal K})$. We
use the notation $\varphi \otimes \eta$ to denote the rank-one operator
$z \mapsto \langle z, \eta \rangle \varphi$. For the purposes of 
computation with Hilbert-Schmidt operators, the formulae
\[ (\varphi \otimes \eta)^* = \eta \otimes \varphi \]
and 
\[ S (\varphi \otimes \eta) T = (S \varphi) \otimes (T^* \eta) \]
are very convenient. As a matter of fact, all calculations in
the algebra ${\cal B}_2({\cal K})$ can be carried out using these
relations, since the rank-one operators span a dense subspace.
Given any projection $P$ and any orthonormal basis $( \eta_i )_{i \in I}$
of $P({\cal K})$, the operators $T \in {\cal B}_2({\cal K}) \circ P$ can 
be shown to have the form 
\begin{equation} \label{eq:dec_T} 
 T = \sum_{i \in I} \varphi_i \otimes \eta_i ~~,
\end{equation}
with $\| T \|_2^2 = \sum_{i \in I} \| \varphi_i \|^2$. Moreover, for
$S =  \sum_{i \in I} \psi_i \otimes \eta_i$, we compute 
\begin{equation} \label{eqn:sc_prod}
\langle S, T \rangle
= \sum_{i \in I} \langle \psi_i, \varphi_i \rangle~~.
\end{equation}
The Hilbert space
$ {\cal B}_2({\cal K}) \circ P$ can be thought of as a tensor product
space or a direct sum of copies of ${\cal K}$:
\[ {\cal B}_2({\cal K}) \circ P \simeq {\cal K} \otimes \ell^2(I) 
 \simeq \left( {\cal K} \right)^{|I|} ~~.\]
This is particularly useful when representations are considered:
If $\pi$ is any representation of the group $G$ on ${\cal K}$,
the representation $\pi \otimes \overline{\pi}$ of $G \times G$ operates
on ${\cal B}_2({\cal K})$ by
\[ (\pi \otimes \overline{\pi})(x,y) T = \pi(x) T \pi(y)^* ~~.\]
Denote by $\pi \otimes 1$ the restriction of this representation to
$G \times \{ 1 \} \simeq G$.
The space ${\cal B}_2({\cal K}) \circ P$ is invariant under
$\pi \otimes 1$, and in the decomposition (\ref{eq:dec_T}) the action
is given as
\[ (\pi \otimes 1) (x) T = \sum_{i \in I} (\pi(x) \varphi_i) \otimes \eta_i ~~.
\] 
This shows that $\pi \otimes 1$, restricted to ${\cal B}_2({\cal K}) \circ P$,
is equivalent to the direct sum of ${\rm rank}(P)$ copies of $\pi$.
Hence ${\rm rank}(P)$ is the multiplicity of $\pi$ in the
restriction of $\pi \otimes 1$ to ${\cal B}_2({\cal K}) \circ P$. 
In particular, given
a set $\Gamma \subset G$, (\ref{eqn:sc_prod}) entails the equivalence
\begin{equation} \label{eqn:equiv_tf}
(\pi \otimes 1) (\Gamma) T \subset {\cal B}_2({\cal K}) \circ P
 \mbox{ is a tight frame} \Longleftrightarrow
 \left( (\pi(\gamma) \varphi_i)_{i \in I} \right)_{\gamma \in \Gamma}
 \subset {\cal K}^{|I|}
 \mbox{ is a tight frame}~~,
\end{equation}
with the same frame constant.
Let us next turn to direct integrals. It is most convenient to define
measurable vector fields first and then measurable operator fields in
terms of the former. Since all Schr\"odinger representations live
on the same space, the direct integral Hilbert space of measurable
vector fields is easily identified with an ordinary ${\rm L}^2$-space:
\[ \int_{\RR'}^{\oplus} {\rm L}^2(\RR) |h| dh \simeq
 {\rm L}^2(\RR,dx) \otimes {\rm L}^2(\RR',|h|dh) \simeq {\rm L}^2(
 \RR \times \RR',dx |h|dh) ~~.\]
Hence measurable vector fields $(\eta^h)_{h \in \RR'}$ may be identified
with families of functions such that $(x,h) \mapsto \eta^h(x)$ is measurable 
and square-integrable with respect to $dx |h|dh$. Measurable operator fields
are operators on ${\rm L}^2(\RR,dx) \otimes {\rm L}^2(\RR',|h|dh)$
of the form $(\eta^h)_h \mapsto (T^h \eta_h)_h$. Now let ${\cal H}
\subset {\rm L}^2(\HH)$ be leftinvariant, with associated
field of projection operators $(\widehat{P}_h)_{h \in \RR}$ and
multiplicity function $m(h) = {\rm rank}(\widehat{P}_h)$.
Define $I_h = \{ 1,2, \ldots, m(h) \}$; by
convention, $I_h = \NN$ for $ m(h) = \infty$.
It is possible to pick a family of vector fields $(\eta_i^h)_{h \in \RR'}$,
with the property that for almost all $h$, $\{\eta^h_1,\ldots, \eta^h_{m(h)}
 \} $ is an orthonormal basis of ${\rm L}^2(\RR)$, and $\eta^h_n = 0$
for $n> m(h)$ \cite[Proposition 7.27]{Fo2}.
Then the elements $f \in {\cal H}$ are characterised on the 
Plancherel transform side by
\[ \widehat{f}(h) = \sum_{i \in I_h} \psi_i^h \otimes \eta_i^h ~~,\]
with measurable vector fields $(\psi_i^h)_{h \in \RR'}$. Hence, in constructing
such elements on the Plancherel transform side, the measurability
of the vector fields $(\psi_i^h)_{h \in \RR'}$ (with $i \in \NN$)
ensures the measurability of the resulting operator field.
As a particular application we note the following lemma:
\begin{lemma} \label{lem:enough_m}
If $m : \RR' \to \NN$ is measurable, there exists a closed
leftinvariant subspace ${\cal H} \subset {\rm L}^2(\HH)$
with $m$ as multiplicity function.
\end{lemma}

\end{appendix}

\section*{Acknowledgements}
The author would like to thank M. Lindner of TU M\"unchen for pointing
out the reference \cite{Ba}, as well as G. Schlichting of TU M\"unchen
for useful discussions.

%File name: bib.tex


\begin{thebibliography}{99}
\bibitem{AlFuKr}{S.T. Ali, H. F\"uhr and A. Krasowska: {\em Plancherel 
inversion as unified approach to wavelet transforms and Wigner functions.}
Submitted to Ann. Inst. H. Poincar\'e. Electronically available as
 \underline{\sf math-ph/0106014}.}
\bibitem{Ba}{R. Balan: {\em Density and redundancy of the noncoherent
Weyl-Heisenberg superframes}, Cont. Mathematics {\bf 247} (1999), 29-41.}
\bibitem{CoGr}{L. Corwin and F.P. Greenleaf: {\em Representations of Nilpotent
 Lie Groups and Their Applications.} Cambridge University Press, Cambridge,
 1989.}
\bibitem{Da} {I. Daubechies: {\em The wavelet transform, time-frequency
 localization and signal analysis.} 
 IEEE Trans. Inform. Theory {\bf 34} (1988), 961-1005.}
\bibitem{Di} {J. Dixmier: {\em $C^{\ast}$-Algebras.}
 North Holland, Amsterdam, 1977.}
\bibitem{Do}{A.H. Dooley: {\em A nonabelian version of the Shannon sampling
theorem.} Siam. J. Math. Anal. {\bf 20} (1989), 624-633.}
\bibitem{Fo} {G.B. Folland: {\em Harmonic Analysis in Phase Space.}
  Princeton University Press, Princeton, 1989.}
\bibitem{Fo2} {G.B. Folland: {\em A Course in Abstract Harmonic Analysis.}
 CRC Press, Boca Raton, 1995.}
%% \bibitem{Fu} {H. F\"uhr: {\em Wavelet frames and admissibility in higher
%% dimensions.} J. Math. Phys. {\bf 37} (1996), 6353-6366.}
\bibitem{FuMa} {H. F\"uhr and M. Mayer: {\em Continuous wavelet transforms
 from semidirect products: Cyclic representations and Plancherel measure.}
 To appear in J. Fourier Anal. Appl. Electronically available as
 \underline{\sf math-ph/0102002}.}
\bibitem{Fu} {H. F\"uhr: {\em Admissible vectors for the regular
 representation.} To appear in Proc AMS. Electronically available as
 \underline{\sf math-ph/0010051}.}
 \bibitem{Gr}{K. Gr\"ochenig, {\em Foundations of Time-Frequency Analysis.}
  Birkh\"auser, Boston, 2001. }
 \bibitem{HiSt}{J.R. Higgins and R.L. Stens: {\em Sampling Theory in Fourier
 and Signal Analysis. Advanced Topics.} Oxford University Press, Oxford, 1999.}
 \bibitem{HeRo} {E. Hewitt and K.A. Ross: {\em Abstract harmonic
  analysis I.} Springer Verlag, Berlin, 1963.}
 \bibitem{Kl}{I. Kluv\'anek, {\em Sampling theorem in abstract
 harmonic analysis.} Mat.-Fyz. Casopis Sloven.Akad. Vied. {\bf 15}
 (1965), 43-48.}
\bibitem{Li} {R.L. Lipsman: {\em Non-abelian Fourier analysis.}
Bull. Sci. Math. {\bf 98} (1974), 209-233.}
%\bibitem{Lu} {J.Ludwig: {\em On the Hilbert-Schmidt Semi-Norms of $L^1$ of a
%nilpotent Lie group}, Math. Ann. {\bf 273} (1986) 383-395.}
%% \bibitem{Li}{R.L. Lipsman: {\em Harmonic analysis on non-semisimple symmetric
%%  spaces.} Israel J. Math. {\bf 54} (1986), 335-350.}
%% \bibitem{LoMaRi} {A.K. Louis, P. Maa\ss{} and A. Rieder: {\em Wavelets.}
%%  Teubner, Stuttgart, 1994.}
%% \bibitem{Maa}{P. Maa\ss{}: {\em Families of two-dimensional wavelets.}
%%  SIAM J. Math. Anal. {\bf 27} (1996), 1454-1481.}
%\bibitem{Ma68} {G. W. Mackey: {\em Induced representations of groups
%   and quantum mechanics.} W. A. Benjamin Inc., New York, 1968.}
%% \bibitem{Ma76} {G. W. Mackey: {\em The theory of unitary group
%% representations.} University of Chicago Press, Chicago, 1976.}
%% \bibitem{Ma} {G. W. Mackey: {\em Point realizations of transformation groups.}
%%  Illinois J. Math. {\bf 6} (1962), 327-335.}
%\bibitem{Ma2}{G. W. Mackey: {\em Induced representations of locally
%  compact groups, I.} Ann. of Math. 
%  {\bf 55} (1952), 101-139.}
%\bibitem{Ma3}{G. W. Mackey: {\em Borel structure in groups and their duals.}
% Trans. Amer. Math. Soc. {\bf 85} (1957), 134-165.}
%% \bibitem{May} {M. Mayer: {\em Square-integrability of tensor products and
%%  symmetry in Fourier-Stieltjes algebras.} Preprint.}
%% \bibitem{McM}{J. McMullen:{\em Extensions of positive definite functions,}
%%  Mem. Amer. Math. Soc. {\bf 117} (1972)}
%% \bibitem{Mu} {R. Murenzi: {\em Ondelettes multidimensionelles et 
%%  application \`{a} l'analyse d'images}. Th\`{e}se, 
%%  Universit\'{e} Catholique de Louvain, Louvain-La-Neuve, 1990.}
%% \bibitem{Ni} {O.A.~Nielsen: {\em Direct integral theory.} Marcel Dekker,
%%  New York, 1980.}
  \bibitem{Ri} {M.~Rieffel: {\em Square integrable representations of 
  Hilbert algebras.} J. Funct. Anal. {\bf 3} (1969), 265--300.}
% \bibitem{Ru} {W. Rudin: {\em Functional Analysis.}
%  McGraw-Hill, New York, 1973.}
% \bibitem{Ru2} {W. Rudin: {\em Real and Complex Analysis.} McGraw-Hill,
% New York, 1966.}
%\bibitem{Wig} {E.P. Wigner: {\em On the quantum correction for thermodynamic 
%  equilibrium.} Phys. Rev. {\bf 40} (1932), 749-759.}
%\bibitem{Wo} {K.B.\ Wolf: {\em Wigner distribution function 
%     for paraxial polychromatic optics.} 
%     Opt.\ Comm. {\bf132} (1996), 343--352.}

%% \bibitem{Zi} {
%%  R. J. Zimmer: {\em Ergodic theory and semisimple groups.}
%%  Birkh\"{a}user, Boston, 1984.}
\end{thebibliography}
\end{document}